\documentclass{article}
\usepackage{amsmath,amsfonts,epsfig}
 \newcommand{\se}{\setcounter{equation}{0}}

\usepackage{amsmath,amsfonts,amssymb}
\usepackage{amsmath,amsfonts,epsfig}
\usepackage{graphics}
\usepackage{enumitem}
\usepackage{supertabular}
\usepackage{amsmath,amsfonts,amssymb}
\usepackage{amsthm}
\usepackage{mathtools}
\usepackage{amsfonts}
\usepackage{cite}
\usepackage{graphicx}
\usepackage{epstopdf}
\usepackage[numbers]{natbib}
\usepackage{float,epsfig, floatflt,here}
\usepackage{color}
\usepackage[colorlinks=true,citecolor=blue,linkcolor=blue]{hyperref}
\usepackage{fancyhdr}
\usepackage{etoolbox}
\usepackage{ dsfont }
\usepackage{accents}
\newcommand{\sKT}{\sum_{K\in\mathcal{T}_h}}

\newcommand{\tnorm}[1]{{\left\vert\kern-0.25ex\left\vert\kern-0.25ex\left\vert #1\right\vert\kern-0.25ex\right\vert\kern-0.25ex\right\vert}}

\newcommand{\vertiii}[1]{{\left\vert\kern-0.25ex\left\vert\kern-0.25ex\left\vert #1
		\right\vert\kern-0.25ex\right\vert\kern-0.25ex\right\vert}}

\newtheorem{them}{Theorem}[subsection]

\newtheorem{lemma}{Lemma}[section]
\newtheorem{lema}{Lemma}[subsection]

\newtheorem{remm}{Remark}[subsection]

\newtheorem{example}{Example}[section]

\begin{document}
\title{\normalsize
  A Systematic Study on Weak Galerkin Finite Element Method for Second Order Parabolic Problems}
\author{{\normalsize Bhupen Deka}\thanks{Department of Mathematics,
			Indian Institute of Technology Guwahati,
			Guwahati - 781039, India (bdeka@iitg.ac.in).}
		{ \normalsize  and  Naresh Kumar}\thanks{Department of Mathematics,
			Indian Institute of Technology Guwahati,
			Guwahati - 781039, India (nares176123101@iitg.ac.in).}}
 \date{}
\maketitle
\begin{abstract}
A systematic numerical study on weak Galerkin (WG) finite element method for second order linear parabolic problems is presented by allowing polynomial approximations with various degrees for each local element.
Convergence of both semidiscrete and fully discrete WG solutions are established
in $L^{\infty}(L^2)$ and $L^{\infty}(H^1)$ norms for a general WG element $({\cal P}_{k}(K),\;{\cal P}_{j}(\partial K),\;\big[{\cal P}_{l}(K)\big]^2)$, where $k\ge 1$, $j\ge 0$ and $l\ge 0$ are arbitrary integers.
The fully discrete space-time discretization is based on a first order in time Euler scheme.
Our  results  are  intended  to  extend  the  numerical analysis of WG methods for elliptic problems [J. Sci. Comput., 74 (2018), 1369-1396] to parabolic problems. Numerical experiments are reported to justify the robustness,
reliability and accuracy of the WG finite element method.
\end{abstract}

\vspace{.01in}


  \section{\normalsize Introduction}\label{s2}\se
 We consider the following linear parabolic equation of the form
\begin{equation}
	u_t-\nabla \cdot(a\nabla u)=f\;\; \mbox{in}\;\Omega \times (0,T]\label{1:1}
\end{equation}
with initial and boundary conditions
\begin{equation}
	u(x,0)=\psi(x)\;\;\mbox{in}\; \Omega;\;\; u=0
	\;\;\;\mbox{on} ~\partial\Omega \times
	(0,T], \label{1:2}
\end{equation}
where $\Omega \subset \mathbb{R}^2$ is a bounded domain with smooth
boundary $\partial\Omega$.
We assume that the coefficient matrix ${a}=(a_{ij}(x))_{2\times 2}\in [L^{\infty}(\Omega)]^{2^2}$
is symmetric and uniformly positive definite in
$\Omega$. The initial function $\psi=\psi(x)$
and the forcing function $f=f(x, t)$ are assumed to be smooth functions in their respective domains of definition, and $T$ is the finite terminal observation time.

Finite element approximations of linear parabolic equations have been studied extensively, \cite{thomee1990finite} contains a comprehensive list of references.
Recently, the weak Galerkin finite element method has attracted much
attention in the field of numerical partial differential equations. The objective of the present work is to propose a systematic framework for the weak Galerkin
finite element method (WG-FEM for short) for second order linear parabolic equations by using polynomials of
various degrees in the weak finite element space. The WG-FEM introduced in \cite{wang2013weak} refers to the numerical algorithms for differential
equations where the differential operators appearing in the variational forms
are to project into another appropriately chosen Sobolev space
such that its approximation by polynomials is possible.
 More precisely, the WG finite element approximations are derived from the weak formulations of the problems by replacing differential operators involved by its weak forms and adding parameter free stabilizers.
 In fact, WG formulation is a natural extension of conforming finite element formulation when nonconforming elements are used. The
concept of weak derivatives makes WG a widely applicable numerical
technique for a large variety of of PDEs arising from the mathematical modeling of practical problems in science and engineering.
There is an abundant literature on such PDEs; see, e.g., elliptic equation
\cite{liu2020penalty, li2018interior, wang2014weak, wang2018systematic, wang2018primal, mu2017least, lin2018weak, liu2018lowest},
parabolic equation \cite{xie2019error, zhang2016weak, li2013weak, zhou2019weak}, system of
equations \cite{mu2020uniformly, mu2020pressure, wang2016discretization, wang2016weak, zhang2019globally, mu2015weak, shields2017weak, liu2016weak, mu2018weak},
interface problems \cite{deka2019weak, song2018relaxed, mu2016new}.
  One close relative of the WG finite element method is the
hybridizable discontinuous Galerkin (HDG) method \cite{cockburn2009unified}. But they make use of different polynomial
approximating spaces and utilize different stabilization techniques. For detailed discussions, we refer to \cite{wang2019basics, cockburn2018weak}.

 The classical finite element methods based on conforming finite element discretization, have limitations in practical computation. The conforming finite element space
is restricted to piecewise polynomials with prescribed continuity that ensures conformity and
stability of the corresponding weak formulation, which is often very difficult to implement, particularly for problems in
high dimensions and/or on general polytopal partitions. In scientific computing, higher order
of convergence is always one of the major research goals, because
high order methods are more accurate and cost efficient. Although conforming finite element
methods have simple formulations with many fewer unknowns,
construction of conforming finite element spaces of any orders
would be either challenging or impossible.
Keeping in
mind the applicability of numerical methods of higher order with polygonal meshes, recently attempts have been made to develop
certain technologies which make use of polygonal meshes, for instance, see \cite{da2014virtual, vacca2015virtual, beirao2013basic} for virtual element methods,
\cite{cangiani2017hp, cangiani2014hp, cockburn2009unified, arnold2002unified, cockburn2016bridging} for discontinuous Galerkin methods.
 Due to the use of discontinuous approximation functions, WG-FEMs are highly flexible in construction of finite element spaces of any orders with the price of more degrees of freedom and more complex formulations.
  Unlike classical finite element method, the WG-FEM is applicable for unstructured polygonal meshes making it
more suitable for complex geometry usually appeared in real life problem.
  A typical local WG element is of the form $({\cal P}_{k}(K),\;{\cal P}_{j}(\partial K),\;\big[{\cal P}_{l}(K)\big]^2)$, where $k\ge 1$ is the degree of polynomials in the interior of the element
  $K$, $j\ge 0$ is the degree of polynomials on the boundary of $K$, and $l\ge 0$ is the degree of polynomials employed in the computation of weak gradients or weak first order partial derivatives.
  The accuracy and the computational complexity of the corresponding WG scheme is significantly impacted by the selection of such polynomials. The goal of this study is to explore all possible combinations of
polynomial functions in the reconstruction of the underlying differential operators. Our results are intended
to extend the weak Galerkin analysis in \cite{wang2018systematic} for elliptic problems to linear parabolic equations with polygonal meshes. More precisely, the analysis presented in this article shows
 that the WG finite element
 solutions approximate the true solutions with an optimal
 order in $L^{\infty}(L^2)$ and $L^{\infty}(H^1)$ norms. The results for parabolic equation are particularly useful because it demonstrate the
 robustness of the WG-FEMs with various combinations of polynomials in the numerical scheme and it fills a gap in literature.
It is worth to note that only $H^1$ norm error estimate is established in
 \cite{wang2018systematic}.
Finally, theoretical convergence results are validated for several combination of the polynomial spaces.

The rest of the paper is organized as follows. In Section 2, we introduce some commonly used notations. Further, we review the definitions of weak
gradient and its discrete analogs in suitable polynomial spaces. Section 3 is devoted to
the optimal order error estimates of semidiscrete WG-FEM algorithm. In Section 4, a backward Euler scheme is described along
with a priori error bounds in $L^{\infty}(H^1)$ and $L^{\infty}(L^2)$ norms. Section 5 focuses on some numerical results that confirm the convergence theory developed in earlier section.
Summary on the new results developed in this paper are presented in Section 6. Finally, in ``Appendix" we present some detailed computational results leading to the reported rate of convergence for some typical combinations.

Throughout the paper, $C$ is a positive generic constant independent of
the mesh parameters $\{h, \tau\}$ and whose value changes with context.

\section{\normalsize Preliminaries and Weak Galerkin Discretization}\se
\subsection{\normalsize Basic Notations} \se
Let us introduce some notations used in this paper. In this work, we will follow the standard notation for Sobolev spaces and
norms (cf. \cite{ra}). For a domain ${\cal K}\subseteq\Omega\subset \mathbb{R}^2$, non-negative integer $m$ and real $p~(1\le p \le \infty)$,
$W^{m,p}({\cal K})$ represents the standard Sobolev space (cf. \cite{ra}). For $p=2$, we use $H^m({\cal K})$ for $W^{m,2}(\cal K)$ with inner product $(\cdot, \cdot)_{m, {\cal K}}$. Notations $\|\cdot \|_{m,{\cal K}}$ and
$|\cdot |_{m, {\cal K}}$ are used to denote the norm and seminorm in the Sobolev space $H^m({\cal K})$, respectively. Inner product in $H^m({\cal K})$ is denoted by
$(\cdot, \cdot)_{m, {\cal K}}$. Clearly, $H^0({\cal K})=L^2({\cal K})$ with inner product $(\cdot, \cdot)_{\cal K}$ and induced norm $\|\cdot \|_{\cal K}$. For our convenience, we skip the subscript ${\cal K}$ in the inner product notation and norm when ${\cal K}=\Omega.$ $H_{0}^{1}(\Omega )$ is the collection of all $H^1(\Omega)$ functions vanishing on the boundary of $\Omega$.

For a given Banach space $(\cal B, \|.\|_{\cal B})$ and interval $J\subset \mathbb{R}$, we define for $m=0,1$,
$$
H^{m}(J;{\cal B})=\Bigg(u(t)\in {\cal B}\;\mbox{for a.e.}\; t\in J
\mbox{ and } \sum_{j=0}^m\int_J\left\|\frac{\partial^j u(t)}{\partial t^j}\right\|^2_{\cal
B}dt<\infty\Bigg)
$$
endowed with the following norm
\begin{eqnarray*}
\|u\|_{H^{m}(J;{\cal B})}=\Bigg(\sum_{j=0}^m\int_J\left\| \frac{\partial^j u
(t)}{\partial t^j}\right\|^2_{\cal
B}dt\Bigg)^{\frac{1}{2}}.
\end{eqnarray*}
Further, $L^{\infty}(J; {\cal B})$ is also a Banach space with respect to following norm
  \begin{eqnarray*}
    \|\phi\|_{L^{\infty}(J; {\cal B})}:= \mbox{ess}\sup_{t\in [0,T]}\|\phi(t)\|_{\cal B}.
  \end{eqnarray*}
For the simplicity, we use $L^{2}({\cal B})$ for $L^{2}(J;{\cal B})$, $L^{\infty}({\cal B})$ for $L^{\infty}(J;{\cal B})$ and $H^{1}({\cal B})$ for $H^{1}(J;{\cal B})$.

We end this section with the following regularity result for the initial boundary value problem (\ref{1:1})-(\ref{1:2}) (see, \cite{necas1996weak}, p. 287, Theorem 2.10).
\begin{them}
Assume that $f\in L^2(0, T;H^{r-1}(\Omega))$ and $\psi \in H^{r}(\Omega)$ for some $r\ge 1$. Then the solution of (\ref{1:1})-(\ref{1:2}) satisfies
$$u\in L^2(0, T;H^{r+1}(\Omega))\cap H^1(0, T;H^{r-1}(\Omega)).\;\;\;\Box$$
\end{them}

\subsection{\normalsize Weak Galerkin Discretization}\se
  In this section, we shall describe the
 weak Galerkin finite element discretization for
  the problem (\ref{1:1})-(\ref{1:2}) and review the definition of the weak gradient operator.

Let ${\cal T}_h$ be a partition of the domain $\Omega$ consisting of polygons in two dimension satisfying a set of conditions specified in \cite{wang2018systematic, wang2014weak}. Denote by ${\cal E}_h$ the set of all edges in ${\cal T}_h$ and let ${\cal E}_h^0={\cal E}_h\backslash \partial \Omega$ be the set of all interior edges. For every element $K\in {\cal T}_h$, we denote by $|K|$ the measure of $K$ and by $h_K$ its diameter and mesh size $h=\max_{K\in {\cal T}_h}h_{K}$ for ${\cal T}_h$.

The key in weak Galerkin methods is the use of weak
derivatives in the place of strong derivatives in the variational form for the underlying
partial differential equations. Thus, it is essential to introduce a weak version for the gradient operator. Weak gradient operators and its discrete version were introduced in \cite{wang2013weak, wang2014weak}, and the rest
of the section will review them. Let $K$ be any polygonal domain with interior $K^0$ and boundary
	$\partial K$.  A weak function on the
region $K$ refers to a pair of scalar valued functions $v=\{v_0, v_b\}$ such that $v_0\in L^2(K)$ and $v_b\in L^2(\partial K)$. 
Denote by ${\cal V}(K)$ the space
	of weak scalar valued functions on $K$; i.e.,
		\begin{equation}
	{\cal V}(K)=\{{v} = \{{v}_0, {v}_b\}: {v}_0 \in L^2(K), {v}_b \in L^2(\partial K)\}.\label{2:1}
	\end{equation}

     For any given integer $k\geq 0$, denote ${\cal P}_k(K)$ the space of polynomials of total degree $k$ or less on the element $K \in {\cal T}_h$. Analogously, for any given integer $j \geq 0$, ${\cal P}_j(e)$ denotes the space of polynomials of total degree $j$ or less on the edge $e \in {\cal E}_h$.
      On each element $K \in {\cal T}_h$, define the following local weak finite element space
   \begin{equation}
	{\cal V}(k,j,K)=\{{v} = \{{v}_0, {v}_b\}: {v}_0 \in {\cal P}_k(K),\; {v}_b \in {\cal P}_j(\partial K)\}. \label{2:2}
	\end{equation}
 A global weak finite element space ${\cal V}_h$ is constructed by patching local space ${\cal V}(k,j,K)$ through a common value of $v_b$ on all interior edges
   \begin{equation}
    {\cal V}_h = \{v=\{v_0, v_b\}:\; v|_K \in {\cal V}(k,j,K),\; [v]_e=0,\;\forall e \in{\cal E}_h^0\}.\label{2:3}
   \end{equation}
Here, $[v]_e$ denotes the jump of $v\in {\cal V}=\prod_{K\in {\cal T}_h}{\cal V}(k,j,K)$ across an interior edge $e\in{\cal E}_h^0$.
     Denote by ${\cal V}_h^0$ the subspace of ${\cal V}_h$ consisting of all finite element functions with vanishing boundary value
     \begin{equation}
     {\cal V}_h^0 = \{v \in {\cal V}_h : v_b|_{\partial \Omega} = 0\}. \label{2:4}
     \end{equation}

      Next, we introduce a discrete weak gradient operator, denoted by $\nabla_{w}$, is defined as the unique
	polynomial $(\nabla_{w}{v}) \in [{\cal P}_{l}(K)]^{2}$ that satisfies the following  equation
	\begin{equation}
	(\nabla_{w}{v},\phi)_K=-\int_{K}{v}_0(\nabla \cdot \phi)dK + \int_{\partial K}{v}_b(\phi \cdot{\bf n})ds\;\;\forall \phi\in[{\cal P}_{l}(K)]^{2}. \label{2:5}
		\end{equation}
     where ${\bf n}$ is the outward normal to $\partial K$ and $l \geq 0$ is prescribed non-negative integer. By applying the divergence theorem to the first term on the right-hand side of (\ref{2:5}) we arrive at
        \begin{equation}
        (\nabla_wv,\phi)_K=(\nabla v_0 ,\phi)_K + \langle v_b-v_0,\phi \cdot {\bf n}\rangle_{\partial K}\;\; \forall \phi \in [P_l(K)]^2. \label{2:6}
       \end{equation}
       Using the discrete weak gradient operator $\nabla_w$, we define a bilinear map ${\cal A}:{\cal V}_h\times  {\cal V}_h \to \mathbb{R}$ by
\begin{equation}	
{\cal A}(u_h,v_h)=\sKT\big(a \nabla_w u_h,\nabla_w v_h\big)_K+{\cal S}(u_h,v_h)\;\;\forall u_h, v_h\in {\cal V}_h.\label{2:6n}
\end{equation}
Here, ${\cal S}(\cdot, \cdot)$ is known as stabilizer, which is a semi-positive definite bilinear form defined on
${\cal V}_h\times {\cal V}_h$. Stabilizer ${\cal S}(\cdot, \cdot)$ is often chosen in such a way that it fits well into the theory and  implementation of the WG numerical scheme. For examples (cf. \cite{wang2018systematic}):
\begin{example} (Projected Element-Boundary Discrepancy) For $v_h=\{v_0, v_b\} \in {\cal V}_h$, the continuity of $v_h$ can be measured by the quantity $v_b-v_0|_{\partial K}$ for each element
    $K \in {\cal T}_h$. The projected element-boundary-discrepancy method is based on
    the following stabilizer
    \begin{equation}
    {\cal S}(u_h,v_h)=\sKT h^{-1}_K \langle {\cal Q}_m(u_b-u_0|_{\partial K}), {\cal Q}_m(v_b-v_0|_{\partial K})\rangle_{\partial K},\label{2:7}
    \end{equation} where  ${\cal Q}_m:L^2(\partial K) \to  {\cal P}_m(\partial K)$
    is the usual $L^2$- projection operator and $m= \max\{j, l\}$.
    \end{example}
 \begin{example} (Element-Boundary Discrepancy) The element-boundary-discrepancy method is based on
    the following stabilizer
    \begin{equation}
    {\cal S}(u_h,v_h)=
    \sKT h^{-1}_K \langle u_b-u_0|_{\partial K}, v_b-v_0|_{\partial K}\rangle_{\partial K}. \label{2:8}
    \end{equation}
\end{example}

In the WG methods, the polynomial degree and the stabilizer must be chosen so that the
bilinear form ${\cal A}(\cdot, \cdot)$ is coercive with respect to the semi-norm $\|\cdot\|_{1, h}$ (cf. \cite{wang2018systematic}) defined by
\begin{equation}
\|v_h\|_{1,h}=\Big(\sum_{K\in {\cal T}_h}(\|\nabla v_0\|_K^2+h_K^{-1}\|v_0-v_b\|_{\partial K}^2)\Big)^{\frac{1}{2}},\;v_h=\{v_0,v_b\}\in {\cal V}_h.\label{2:9}
\end{equation}
More precisely,
there exist constants $C_1, \;C_2>0$ such that for any $v_h\in {\cal V}_h$, the following inequality holds true
\begin{equation}
C_1\|v_h\|_{1, h}^2\le {\cal A}(v_h, v_h) \le C_2\|v_h\|_{1, h}^2. \label{2:10}
\end{equation}

The coercivity inequality (\ref{2:10}), for both the stabilizers on  weak Galerkin space $({\cal P}_{k}(K),\;{\cal P}_{j}(\partial K),\;\big[{\cal P}_{l}(K)\big]^2)$, is stated below (cf. \cite{wang2018systematic}).
\begin{lemma}
Assume that $l\ge k-1$ and $m=\max\{j, l\}$. Then the coercivity inequality (\ref{2:10}) holds true.
\end{lemma}

We end this section with some standard $L^2$ projections.
For each element $K\in {\cal T}_h$ and edge $e\in {\cal E}_h$, operators ${\cal Q}_k^0: L^2(K)\to {\cal P}_k(K)$ and ${\cal Q}_j^b: L^2(e)\to {\cal P}_{j}(e)$
are the usual $L^2$ projections. Denote by ${\cal Q}_h$ the $L^2$ projection onto the finite element space ${\cal V}_h$ such that ${\cal Q}_h|_{K}=\{{\cal Q}_k^0, {\cal Q}_j^b\}$.
In addition to ${\cal Q}_h$, let $\mathbb{Q}_l:[L^2(K)]^2 \to [{\cal P}_{l}(K)]^2$ be an another local $L^2$ projection.

\section{\normalsize Error analysis for the semidiscrete scheme}\se
 This section deals with the error analysis for the spatially discrete
scheme. Optimal order of convergence in both $L^{\infty}(L^2)$ and $L^{\infty}(H^1)$ norms are
established.

A time-dependent weak function $v_h:[0,T]\to {\cal V}_h$ is written as $v_h(t):=\{v_0(t), v_b(t)\}$ and subsequently we define $v_{ht}(t):=\{v_0^{\prime}(t),v_b^{\prime}(t)\}$, where `$\prime$' denotes the time derivatives.
For simplicity, we use $v_h=\{v_0,v_b\}$ for $v_h(t)$ and $v_{ht}=\{v_0^{\prime}, v_b^{\prime}\}$ for $v_{ht}(t)$.

The continuous-time weak Galerkin finite element approximation to
(\ref{1:1})-(\ref{1:2}) can be obtained by seeking $u_h=\{u_0,u_b\}:[0,T]\to {\cal V}_h^0$ satisfying following equation
\begin{eqnarray}
(u_{ht},v_0)+{\cal A}(u_h,v_h)=(f,v_0)\;\;\forall v_h\in {\cal V}_h^0,\label{3:1}
\end{eqnarray}
where $u_h(0)\in {\cal V}_h^0$ is a suitable approximation of the initial function $\psi$.
Well-posedness of the scheme (\ref{3:1}) can be verified from the fact that weak finite element space ${\cal V}_h^0$ is a normed linear space with respect to the triple norm $\tnorm{\cdot}$ defined as
$$\tnorm{v_h}=\sqrt{{\cal A}(v_h, v_h)},\;v_h\in {\cal V}_h^0.$$

As a standard procedure in finite element method, we split our error into two components using an intermediate operator.
We write
	\begin{align*}
	{u}-{u}_h&=({ u} - {\cal Q}_h{ u}) + ({\cal Q}_h{u}-{u}_h).
	\end{align*}
For simplicity, we introduce the following notation
\begin{equation}	
{e}_h(t):=\{{ e}_0(t),{ e}_b(t)\}={ u}_h(t) - {\cal Q}_h{u}(t),\;\;t\in [0, T].\label{3:2}
\end{equation}
Then $e_h$ satisfies following error equation which is crucial for our later analysis.
\begin{lemma} \label{lm:3:1}
Let $e_h$ be the error as defined in (\ref{3:2}). Then, for all $v_h=\{v_0, v_b\}\in {\cal V}_h^0$, we have
\begin{equation}
(e_{ht},v_0)+{\cal A}(e_h,v_h)=l_1(u,v_h)+l_2(u,v_h)+l_3(u,v_h)+{\cal S}({\cal Q}_hu,v_h),\label{main:e}
\end{equation}
where bilinear forms $l_1(\cdot, \cdot),\;l_2(\cdot, \cdot)$ and $l_3(\cdot, \cdot)$ are given by
\begin{eqnarray*}
 		l_1(u,v_h)&=& \sum_{K\in {\cal T}_h}\Big({\mathbb Q}_l(a{\mathbb Q}_l\nabla {\cal Q}_k^0u)-a\nabla u, \nabla v_0\Big)_{K},\\
 		  l_2(u,v_h)&=& \sum_{K\in {\cal T}_h}\langle ({\mathbb Q}_l(a{\mathbb Q}_l\nabla {\cal Q}_k^0u)- a\nabla u)\cdot {\bf n}, v_b-v_0\rangle_{\partial K},\\
 		 l_3(u,v_h)&=& \sum_{K\in {\cal T}_h}\langle {\cal Q}_j^bu-{\cal Q}_k^0u, {\mathbb Q}_l(a\nabla_wv)\cdot {\bf n}\rangle_{\partial K}.
 	\end{eqnarray*}
\end{lemma}

\noindent \textit{Proof.} For any $v_h=\{v_0, v_b\}\in {\cal V}_h^0$, we test equation (\ref{1:1}) against $v_0$ on each element $K\in {\cal T}_h$ to obtain
\begin{eqnarray}
(f,v_0)&=&(u_t,v_0)-\sum_{K\in {\cal T}_h}(\nabla\cdot(a \nabla u),v_0)_K\nonumber\\
&=& ({\cal Q}_hu_t,v_0)+\sum_{K\in {\cal T}_h}(a \nabla u,\nabla v_0)_K-\sum_{K\in {\cal T}_h}\langle a \nabla u \cdot {\bf n}, v_0\rangle_{\partial K}\nonumber\\
&=& (({\cal Q}_hu)_t,v_0)+\sum_{K\in {\cal T}_h}(a \nabla u,\nabla v_0)_K-\sum_{K\in {\cal T}_h}\langle a \nabla u \cdot {\bf n}, v_0-v_b\rangle_{\partial K},~~~~~~~\label{3:4}
\end{eqnarray}
where we have used the divergence theorem and the fact that $$\sum_{K\in {\cal T}_h}\langle a \nabla u \cdot {\bf n},v_b\rangle_{\partial K}=0.$$
Combining (\ref{3:1}) and (\ref{3:4}), we have
\begin{eqnarray}
(u_{ht},v_0)+{\cal A}(u_h,v_h)&=& (({\cal Q}_hu)_t,v_0)+\sum_{K\in {\cal T}_h}(a \nabla u,\nabla v_0)_K\nonumber\\
&& -\sum_{K\in {\cal T}_h}\langle a \nabla u \cdot {\bf n}, v_0-v_b\rangle_{\partial K}.\label{3:5}
\end{eqnarray}
Then integration by parts together with the identity (\ref{2:6}) and the definition of $\mathbb{Q}_l$ operator yields
  \begin{eqnarray*}
   &&(a\nabla_w{\cal Q}_hu,\nabla_wv)_K\\&&=(\nabla_w{\cal Q}_hu, {\mathbb Q}_l(a\nabla_wv))_K\nonumber\\
    &&= \Big(\nabla {\cal Q}_k^0u, {\mathbb Q}_l(a\nabla_{w}v)\Big)_K + \langle {\cal Q}_j^bu-{\cal Q}_k^0u, {\mathbb Q}_l(a\nabla_{w}v)\cdot {\bf n}\rangle_{\partial K}\nonumber\\
    &&= \Big({\mathbb Q}_l(a {\mathbb Q}_l(\nabla {\cal Q}_k^0u)), \nabla_{w}v)\Big)_K+\langle {\cal Q}_j^bu-{\cal Q}_k^0u, {\mathbb Q}_l(a\nabla_{w}v)\cdot {\bf n}\rangle_{\partial K}\nonumber \\
      &&= \Big({\mathbb Q}_l(a{\mathbb Q}_l(\nabla Q_k^0u)),\nabla v_0)\Big)_K\nonumber + \langle Q_j^bu-Q_k^0u,{\mathbb Q}_l(a\nabla_{w}v)\cdot {\bf n}\rangle_{\partial K} \\&&
    \;\;\;\;+ \langle v_b-v_0,{\mathbb Q}_l(a{\mathbb Q}_l(\nabla Q_k^0u))\cdot {\bf n}\rangle_{\partial K},
  \end{eqnarray*}
so that
\begin{eqnarray*}	
{\cal A}({\cal Q}_hu,v_h)&=&\sKT\big(a \nabla_w {\cal Q}_hu,\nabla_w v_h\big)_K+{\cal S}({\cal Q}_hu,v_h)\nonumber\\
&=&\sKT\Big({\mathbb Q}_l(a{\mathbb Q}_l(\nabla Q_k^0u)),\nabla v_0)\Big)_K + \sKT\langle Q_j^bu-Q_k^0u,{\mathbb Q}_l(a\nabla_{w}v)\cdot {\bf n}\rangle_{\partial K} \nonumber\\&&
    + \sKT\langle v_b-v_0,{\mathbb Q}_l(a{\mathbb Q}_l(\nabla Q_k^0u))\cdot {\bf n}\rangle_{\partial K}+{\cal S}({\cal Q}_hu,v_h),
\end{eqnarray*}
and hence,
\begin{eqnarray}
&&(({\cal Q}_hu)_t,v_0)+{\cal A}({\cal Q}_hu,v_h)\nonumber\\
&&=(({\cal Q}_hu)_t,v_0)+\sKT\Big({\mathbb Q}_l(a{\mathbb Q}_l(\nabla Q_k^0u)),\nabla v_0)\Big)_K
\nonumber\\&&\;\;\;\;+ \sKT\langle Q_j^bu-Q_k^0u,{\mathbb Q}_l(a\nabla_{w}v)\cdot {\bf n}\rangle_{\partial K} \nonumber\\&&
    \;\;\;\;+ \sKT\langle v_b-v_0,{\mathbb Q}_l(a{\mathbb Q}_l(\nabla Q_k^0u))\cdot {\bf n}\rangle_{\partial K}+{\cal S}({\cal Q}_hu,v_h). \label{3:6}
\end{eqnarray}
Subtracting (\ref{3:5}) from (\ref{3:6}) leads to desire result. \;\;\;$\Box$

Next, for a shape regular weak Galerkin discretization ${\cal T}_h$, we recall following crucial estimates for the bilinear maps $l_1,\; l_2$ and $l_3$ from literature \cite{wang2018systematic}.
 \begin{lemma}\label{lm:3:2}
  		Let $\sigma = min\{l+1,k\}$. Assume that $u\in H^{\sigma+1}(\Omega)\cap H^1_0(\Omega)$ then the following estimate holds true
  		\begin{eqnarray*}
  		|l_1(u, v_h)| \leq Ch^\sigma\|u\|_{\sigma+1}\|v_h\|_{1,h}\;\;\forall v_h\in {\cal V}_h^0,
  		\end{eqnarray*}
  		where $C$ is a positive constant depending on $\|a\|_{l+1,\infty}$-the element wise $W^{l+1,\infty}$ norm of the coefficient matrix $a$.
  	\end{lemma}

   \begin{lemma}\label{lm:3:3}
   	 Under the assumptions of Lemma \ref{lm:3:2}, for all $v_h\in {\cal V}_h^0$, we have
   	 \begin{eqnarray*}
   	     |l_2(u, v_h)|\leq Ch^\sigma\|u\|_{\sigma+1}\|v_h\|_{1,h}.
   	 \end{eqnarray*}
   \end{lemma}

  \begin{lemma}\label{lm:3:4}
  	 Let $k, j, l$ be the non-negative integers that define the weak finite element space ${\cal V}_h$. Set $s=min\{k, j\}$ and assume that $s\geq1$. In addition,
  assume that $u\in H^{s+1}(\Omega)\cap H_0^1(\Omega)$ then the following estimate holds true
  	 \begin{eqnarray*}
  	  |l_3(u, v_h)|\leq Ch^s\|u\|_{s+1}\|v_h\|_{1,h}\;\;\forall v_h\in {\cal V}_h^0.
  	 \end{eqnarray*}
  	 In the case $j\geq l$, we have
  	 \begin{eqnarray}
  	    |l_3(u, v_h)|\leq Ch^k\|u\|_{k+1}\|v_h\|_{1,h}\;\;\forall v_h\in {\cal V}_h^0.
  	 \end{eqnarray}
  \end{lemma}

 \subsection{\normalsize Error estimates with  projected element-boundary-discrepancy}
 Convergence results for the semidiscrete weak Galerkin approximation based on projected element-boundary-discrepancy
 are presented.

For our convenience, following result is borrowed from \cite{wang2018systematic}.
 \begin{lema} \label{lm:3:1:1}Assume that
 \begin{equation*}
    {\cal S}(u_h,v_h)=\sKT h^{-1}_K \langle {\cal Q}_m(u_b-u_0|_{\partial K}), {\cal Q}_m(v_b-v_0|_{\partial K})\rangle_{\partial K},
    \end{equation*}
    where ${\cal Q}_m:L^2(\partial K) \to  {\cal P}_m(\partial K)$
    is the usual $L^2$- projection operator and $m= \max\{j, l\}$. The following results hold true:
    \begin{description}
    \item {\rm (a)}\; Assume that the solution of (\ref{1:1})-(\ref{1:2}) is so regular that $u\in H^{k+1}(\Omega)\cap H_0^1(\Omega)$ and $j\ge l,$ then
    $$|{\cal S}({\cal Q}_hu,{\cal Q}_hu)|\le Ch^{2k}\|u\|_{k+1}^2.$$
           \item {\rm (b)}\;Assume that the solution of (\ref{1:1})-(\ref{1:2}) is so regular that $u\in H^{s+1}(\Omega)\cap H_0^1(\Omega)$ and $j<l,$ then
    $$|{\cal S}({\cal Q}_hu,{\cal Q}_hu)|\le Ch^{2s}\|u\|_{s+1}^2,$$ where $s=\min\{k, j\}.$
    \end{description}
 \end{lema}

 The convergence results for the stabilizer with projected element-boundary-discrepancy can be summarized as follows.
 \begin{them}
 Let $k, j,$ and $l\ge {k-1}$ be the non-negative integers that define the weak finite element space ${\cal V}_h$. Assume that
 $${\cal S}(u_h,v_h)=\sKT h^{-1}_K \langle {\cal Q}_m(u_b-u_0|_{\partial K}), {\cal Q}_m(v_b-v_0|_{\partial K})\rangle_{\partial K},$$
  where ${\cal Q}_m:L^2(\partial K) \to  {\cal P}_m(\partial K)$ is the usual $L^2$ projection operator and $m= \max\{j, l\}$. Then the following error estimates hold true:
  \begin{description}
  \item {\rm (a)}\;For $j<l$, set $s=\min\{k, j\}$ and assume $s\geq 1$. Assume that the solution of (\ref{1:1})-(\ref{1:2}) is so regular that $u\in H^{s+1}(\Omega)\cap H_0^1(\Omega)$. Then
    \begin{eqnarray}
 &&\|e_h(t)\|^2+\int_{0}^{t}\|e_h\|^2_{1,h}ds\le C\Big(\|e_h(0)\|^2+h^{2s}\int_{0}^{t}\|u\|_{s+1}^2ds\Big), \label{thm1:a:1}~~~~\\
 &&\int_0^t\|e_{ht}(t)\|^2ds+\|e_{h}(t)\|_{1, h}^2\nonumber\\&&\le C\Big(\tnorm{e_{h}(0)}^2+\|e_{ht}(0)\|^2+ h^{2s}\int_0^t\|u\|_{s+1}^2ds\Big). \label{thm1:a:2}
    \end{eqnarray}
    \item {\rm (b)}\;For $j\ge l$, assume that the solution of (\ref{1:1})-(\ref{1:2}) is so regular that $u\in H^{k+1}(\Omega)\cap H_0^1(\Omega)$. Then
    \begin{eqnarray}
   &&\; \|e_h(t)\|^2+\int_{0}^{t}\|e_h\|^2_{1,h}ds\le C\Big(\|e_h(0)\|^2+h^{2k}\int_{0}^{t}\|u\|_{k+1}^2ds\Big), \label{thm1:b:1}~~~~\\
   &&\; \int_0^t\|e_{ht}(t)\|^2ds+\|e_{h}(t)\|_{1, h}^2\nonumber\\
   &&~~\le C\Big(\tnorm{e_{h}(0)}^2+\|e_{ht}(0)\|^2+ h^{2k}\int_0^t\|u\|_{k+1}^2ds\Big). \label{thm1:b:2}
    \end{eqnarray}
  \end{description}
\end{them}

\noindent \textit{Proof.} Set $v_h=e_h$ in the error equation (\ref{main:e}) to obtain
\begin{eqnarray*}
\frac{1}{2}\frac{d}{dt}\|e_h(t)\|^2+{\cal A}(e_h, e_h)\le |l_1(u,e_h)|+|l_2(u,e_h)|+|l_3(u,e_h)|+|{\cal S}({\cal Q}_hu,e_h)|.
\end{eqnarray*}
Then by integrating from $0$ to $t$ and using the coercive inequality (\ref{2:10}), we have
\begin{eqnarray}
\frac{1}{2}\|e_h(t)\|^2+C_1\int_0^t\|e_h\|_{1, h}^2ds&\le& \int_0^t|l_1(u,e_h)|ds+\int_0^t|l_2(u,e_h)|ds\nonumber\\
&&+\int_0^t|l_3(u,e_h)|ds+\int_0^t|{\cal S}({\cal Q}_hu,e_h)|ds\nonumber\\
&:=& I_1+I_2+I_3+I_4.\label{3:8}
\end{eqnarray}
For the terms $I_1$ and $I_2$, we first observe that $\sigma=\min\{l+1, k\}=k$. Then, we apply Lemma \ref{lm:3:2} and Lemma \ref{lm:3:3} to have
\begin{eqnarray}
I_1,\;I_2\le Ch^k\int_0^t\|u\|_{k+1}\|e_h\|_{1, h}ds. \label{3:9}
\end{eqnarray}
Assume that $j<l$ and $s=\min\{k, j\},$ so that  Lemma \ref{lm:3:4} and Lemma \ref{lm:3:1:1} yields
\begin{eqnarray}
I_3,\;I_4\le Ch^s\int_0^t\|u\|_{s+1}\|e_h\|_{1, h}ds. \label{3:10}
\end{eqnarray}
Combining (\ref{3:8})-(\ref{3:10}), we have following $L^{\infty}(L^2)$ norm and $L^2(\tnorm{\cdot})$ norm error estimates
\begin{eqnarray}
    \|e_h(t)\|^2+\int_{0}^{t}\|e_h\|_{1, h}^2ds\leq C\Big\{\|e_h(0)\|^2+h^{2s}\int_{0}^{t}\|u\|_{s+1}^2ds\Big\}. \label{3:11}
    \end{eqnarray}
    In the last inequality, we have used the standard Young's inequality.

Next, we differentiate (\ref{main:e}) with respect to $t$ and then set $v_h=e_{ht}$ in the resulting equation to have
\begin{eqnarray*}
\frac{1}{2}\frac{d}{dt}\|e_{ht}(t)\|^2+{\cal A}(e_{ht}, e_{ht})&\le& |l_1(u,e_{ht})|+|l_2(u,e_{ht})|\nonumber\\&&
+|l_3(u,e_{ht})|+|{\cal S}({\cal Q}_hu,e_{ht})|.
\end{eqnarray*}
Arguing as in (\ref{3:11}), we note that
\begin{eqnarray}
\|e_{ht}(t)\|^2+\int_0^t\|{e_{ht}}\|_{1, h}^2ds\le C\Big(\|e_{ht}(0)\|^2+ h^{2s}\int_0^t\|u\|_{s+1}^2ds\Big).\label{n:s:1}
\end{eqnarray}

Next, we set $v_h=e_{ht}$ in the error equation (\ref{main:e}) and arguing as in (\ref{3:11}), we obtain following error estimate
\begin{eqnarray*}
\int_0^t\|e_{ht}(t)\|^2ds+\|e_{h}(t)\|_{1, h}^2\le C\Big(\tnorm{e_{h}(0)}^2+\|e_{ht}(0)\|^2+ h^{2s}\int_0^t\|u\|_{s+1}^2ds\Big).
\end{eqnarray*}
Here, we have used the estimate (\ref{n:s:1}).

Part (b) can be realized in a similar manner. We omit the details. This completes the rest of the proof. \;\;\;$\Box$

Next, we derive an optimal order of estimate for $e_h$ in $L^2$ norm, the basic idea applied is to use elliptic projection.
For $v\in H^2(\Omega)\cap H_0^1(\Omega)\}$, we define
\begin{eqnarray*}
f_v=-\nabla\cdot(a\nabla v)\;\;\mathrm{in}\;\Omega.
\end{eqnarray*}
Clearly, $f_v\in L^2(\Omega)$. Define ${\cal R}_h: H^2(\Omega)\cap H_0^1(\Omega)\to {\cal V}_h^0$ by
\begin{eqnarray}
{\cal A}({\cal R}_hv,v_h)=(f_v,v_h)\;\;\forall v_h=\{v_0, v_b\}\in {\cal V}_h^0,\;v\in H^2(\Omega)\cap H_0^1(\Omega).~~~\label{ritz:1}
\end{eqnarray}
It is easy to observe from the definition of elliptic projection and equation (\ref{3:1}) that
\begin{equation}
(u_{ht},v_h)+{\cal A}(u_h-{\cal R}_hu,v_h)=(f,v_h)+(\nabla\cdot(a \nabla u),v_h)=(u_t,v_h), \label{ref:1}
\end{equation}
for all $v_h=\{v_0, v_b\}\in {\cal V}_h^0.$ Here, we have used equation (\ref{1:1}).

\begin{remm}
From the identity (\ref{ref:1}), for $u_h(0)={\cal R}_h\psi$, it is easy to see that
\begin{eqnarray*}
(e_{ht}(0), v_h)&=&(u_{ht}(0)-{\cal Q}_h{u_t}(0),v_h)\\&=&(u_t(0)-{\cal Q}_h{u_t}(0),v_h)\;\;\forall v_h=\{v_0, v_b\}\in {\cal V}_h^0,
\end{eqnarray*}
which implies
\begin{eqnarray}
\|e_{ht}(0)\|\le \|u_t(0)-{\cal Q}_h{u_t}(0)\|\le Ch^{\lambda}\|u_t(0)\|_{\lambda},\;\;0\le \lambda \le k.\label{3:19:nd}
\end{eqnarray}
Here, we have used standard approximation properties for $L^2$ projection (see, Lemma 4.1 in \cite{wang2014weak}). Again, from the equation (\ref{1:1}) it follows that
\begin{eqnarray}
\|u_t(0)\|_{\lambda}\le C\big(\|\psi\|_{\lambda+2}+\|f\|_{H^1(J;H^{\lambda})}\big).\label{3:20:d}
\end{eqnarray}
Combining estimates (\ref{3:19:nd}) and (\ref{3:20:d}), we obtain
\begin{eqnarray}
\|e_{ht}(0)\|\le Ch^{\lambda}\big(\|\psi\|_{\lambda+2}+\|f\|_{H^1(J;H^{\lambda})}\big),\;\;0\le \lambda \le k.\label{3:19:d}
\end{eqnarray}
\end{remm}

In view of (\ref{ritz:1}), we observe that ${\cal R}_hv$ is the WG finite element approximation of the elliptic
problem with exact solution $v\in H^2(\Omega)\cap H_0^1(\Omega)$ satisfying following equation
\begin{eqnarray}
-\nabla \cdot(a \nabla v)=f_v\;\;\; \mbox{in}\;\Omega.\label{2:e:1}
\end{eqnarray}
Then the error $\rho_v:={\cal Q}_hv-{\cal R}_hv$ satisfies following error equation (see, Lemma 4.1 in \cite{wang2018systematic})
\begin{equation}
{\cal A}(\rho_v,w_h)=l_1(v,w_h)+l_2(v,w_h)+l_3(v,w_h)+{\cal S}({\cal Q}_hv,w_h),\label{main:ee}
\end{equation}
for all $w_h\in {\cal V}_h^0$.

Further, following discrete $H^1$ norm error estimates for $R_h$ hold true \cite{wang2018systematic}.
\begin{lema}
 Let $k, j,$ and $l\ge {k-1}$ be the non-negative integers that define the weak finite element space ${\cal V}_h$. Assume that
 $${\cal S}(u_h,v_h)=\sKT h^{-1}_K \langle {\cal Q}_m(u_b-u_0|_{\partial K}), {\cal Q}_m(v_b-v_0|_{\partial K})\rangle_{\partial K},$$
  where ${\cal Q}_m:L^2(\partial K) \to  {\cal P}_m(\partial K)$ is the usual $L^2$ projection operator and $m= \max\{j, l\}$. Then the following error estimates hold true:
  \begin{description}
  \item {\rm (a)}\;For $j<l$, set $s=\min\{k, j\}$ and assume $s\geq 1$. For $v\in H^{s+1}(\Omega)\cap H_0^1(\Omega)$, we have
    \begin{eqnarray}
    \|{\cal Q}_hv-{\cal R}_hv\|_{1, h}\le Ch^s\|v\|_{s+1}. \label{lem1:a:1}
    \end{eqnarray}
    \item {\rm (b)}\;For $j\ge l$ and $v\in H^{k+1}(\Omega)\cap H_0^1(\Omega)$. Then
    \begin{eqnarray}
    \|{\cal Q}_hv-{\cal R}_hv\|_{1, h}\le Ch^k\|v\|_{k+1}. \label{lem1:b:2}
    \end{eqnarray}
  \end{description}
\end{lema}

 Next, the error $e_h=u_h-{\cal Q}_hu$ is expressed in terms of standard $\rho$ and $\theta$ as
\begin{equation}
e_h(t)=u_h(t)-{\cal Q}_hu(t)=\theta(t)-\rho(t), \label{rho:1}
\end{equation}
where $\rho:={\cal Q}_hu-{\cal R}_hu$ and $\theta:=u_h-{\cal R}_hu$.

For $\theta\in {\cal V}_h^0$, we note that (cf. \cite{deka2019weak})
\begin{eqnarray}
(\theta_t,v_h)+{\cal A}(\theta,v_h)=(\rho_t,v_h)\;\;\forall v_h\in {\cal V}_h^0.\label{theta:1}
\end{eqnarray}
For $v_h=\theta$ in (\ref{theta:1}), we have
\begin{eqnarray*}
(\theta_t,\theta)+\tnorm{\theta}^2
\le   \|\rho_t\|\|\theta\|,
\end{eqnarray*}
which leads to
\begin{equation*}
\|\theta\|^2+\int_0^t\tnorm{\theta}^2ds\le \|\theta(0)\|^2+ C\int_0^t\|\rho_t\|^2ds+C\int_0^t\|\theta\|^2ds.
\end{equation*}
A simple application of Grownwall's inequality yields
\begin{equation}
\|\theta\|^2\le C\Big(\|\theta(0)\|^2+\int_0^t\|\rho_t\|^2ds\Big)=C\int_0^t\|\rho_t\|^2ds,\label{theta:2}
\end{equation}
where we have used the fact that $\theta(0)=u_h(0)-{\cal R}_hu(0)=0$.

\begin{remm}
To the best of our knowledge, optimal error estimates in $L^2$ norm for elliptic problems on general WG finite element space $$({\cal P}_{k}(K),\;{\cal P}_{j}(\partial K),\;\big[{\cal P}_{l}(K)\big]^2)$$
 with arbitrary non-negative integers $\{k,\;j,\;l\}$ have not been established
 earlier.
Article \cite{wang2018systematic} is only concerned about the discrete $H^1$ norm convergence. Therefore, we directly can not use optimal convergence results for the term $\rho_t$ in the $L^2$ norm.
\end{remm}

Next, for the $L^2$ norm error estimate, we now consider the following
auxiliary problem: For every $t\in [0, T]$, find $z(t) \in H_0^1(\Omega)\cap H^2(\Omega)$ such that
\begin{equation}
-\nabla \cdot (a \nabla z(t))=\rho_t(t).\label{3:1:9:tc}
\end{equation}
 Then, we may define $z_h(t):=\{z_0(t), z_b(t)\}\in {\cal V}_h^0$ as the solution to following discrete elliptic problem
\begin{equation}
{\cal A}(z_h(t), v_h) = (\rho_t(t), v_h)\;\;\forall v_h\in
V_h^0,\;\;t\in [0, T].\label{3:1:10:tc}
\end{equation}
Clearly, $z_h$ is the weak Galerkin finite element approximation to $z$ and satisfies following estimates (cf. \cite{wang2018systematic})
\begin{equation}
\|{z-z_h}\|_{1, h}\le Ch\|z\|_2\le Ch\|\rho_t\|. \label{jo:1}
\end{equation}
 Here, we have used the standard a priori estimate for elliptic problem and the WG space $({\cal P}_{k}(K),\;{\cal P}_{j}(\partial K),\;\big[{\cal P}_{l}(K)\big]^2)$ with $k$, $j$ and $l\ge k-1$ are non-negative integers.

Setting $v_h=\rho_t$ in (\ref{3:1:10:tc}) and further using
identity (\ref{main:ee}), we have
\begin{eqnarray}
\|\rho_t\|^2 &=& {\cal A}(z_h, \rho_t) \nonumber \\ &=& l_1(u_t, z_h)+l_2(u_t, z_h)+l_3(u_t, z_h)+{\cal S}(Q_hu_t, z_h).\label{3:1:11:tcc}
\end{eqnarray}
Hence, integrating $(\ref{3:1:11:tcc})$ from $0$ to
$T$, we arrive at following estimate
\begin{eqnarray}
\int_0^T\|\rho_t\|^2ds&\le& \int_0^Tl_1(u_t, z_h)ds+\int_0^Tl_2(u_t, z_h)ds\nonumber\\&&+\int_0^Tl_3(u_t, z_h)ds
+\int_0^T{\cal S}({\cal Q}_hu_t, z_h)ds\nonumber\\
&:=& I_1+I_2+I_3+I_4. \label{pa:1}
\end{eqnarray}
We now estimate each term separately. For the term $I_1$, we use the definition of $L^2$ projection and the fact that $\nabla z_0\in [{\cal P}_{k-1}(K)]^2 \subseteq [{\cal P}_l(K)]^2$ to have
\begin{eqnarray}
 		|l_1(u_t,z_h)|&=& \Bigg{|}\sum_{K\in {\cal T}_h}\Big({\mathbb Q}_l(a{\mathbb Q}_l\nabla {\cal Q}_k^0u_t)-a\nabla u_t, \nabla z_0\Big)_{K}\Bigg{|}\nonumber\\
                    &\le& \sum_{K\in {\cal T}_h}\Bigg{|}\Big(a{\mathbb Q}_l\nabla {\cal Q}_k^0u_t-a\nabla u_t, \nabla z_0\Big)_{K}\Bigg{|}\nonumber\\
                    &\le& \sum_{K\in {\cal T}_h}\Bigg{|}\Big(a{\mathbb Q}_l\nabla {\cal Q}_k^0u_t-a\nabla{\cal Q}_k^0u_t, \nabla z_0\Big)_{K}\Bigg{|}\nonumber\\
                    &&+ \sum_{K\in {\cal T}_h}\Bigg{|}\Big(a\nabla {\cal Q}_k^0u_t-a\nabla u_t, \nabla z_0\Big)_{K}\Bigg{|}\nonumber\\
                    &:=&I_{11}+I_{12}. \label{k:1}
 		  \end{eqnarray}
 Now, we use approximation properties for $L^2$ projections to have
 \begin{eqnarray}
 		I_{11}&=&\sum_{K\in {\cal T}_h}\Bigg{|}\Big({\mathbb Q}_l\nabla {\cal Q}_k^0u_t-\nabla{\cal Q}_k^0u_t, (a -{\bar a})\nabla z_0\Big)_{K}\Bigg{|}\nonumber\\
 &\le& Ch\|a\|_{1, \infty}\sum_{K\in {\cal T}_h}\Bigg{|}\Big({\mathbb Q}_l\nabla {\cal Q}_k^0u_t-\nabla{\cal Q}_k^0u_t, \nabla z_0\Big)_{K}\Bigg{|}\nonumber\\
 &\le& Ch\|a\|_{1, \infty}\sum_{K\in {\cal T}_h}Ch^{\lambda_1+1}\|\nabla {\cal Q}_k^0u_t\|_{\lambda_1+1, K}\|\nabla z_0\|_{K}\nonumber\\
  &\le& Ch^{\lambda_1+2}\|a\|_{1, \infty}\|u_t\|_{\lambda_1+2}\|z_h\|_{1, h}.\label{I:1:1}
		  \end{eqnarray}
 for some non-negative integer $\lambda_1\le l$. Here, ${\bar a}$ is the average of $a$ on each element $K\in {\cal T}_h$.

 For the term $I_{1 2}$, we use the shape regularity assumptions described in \cite{wang2014weak}. For any $K\in {\cal T}_h$, we have a shape regular circumscribed simplex  $S(K)$
 with diameter $h_{S(K)}$ such that $h_{S(K)}\le \gamma_\star h_K$ with a constant $\gamma_{\star}>0$ independent of $K\in{\cal T}_h$. The shape regularity of $S(K)$ implies
 that the measure of $S(K)$ is proportional to $h_{S(K)}^2$.
  Then, we obtain
 \begin{eqnarray}
 		I_{12}&\le& \sum_{K\in {\cal T}_h}\|a\|_{L^{\infty}(K)}\|\nabla {\cal Q}_k^0u_t-\nabla u_t\|_K
                    \|\nabla z_0\|_{K}\nonumber\\
                    &\le& C\sum_{K\in {\cal T}_h}\|a\|_{2, K}\|\nabla {\cal Q}_k^0u_t-\nabla u_t\|_K
                    \|\nabla z_0\|_{K}\nonumber\\
                    &\le& C\sum_{K\in {\cal T}_h}|K|^{\frac{1}{2}}\|a\|_{W^{2, \infty}(K)}\|\nabla {\cal Q}_k^0u_t-\nabla u_t\|_K
                    \|\nabla z_0\|_{K}\nonumber\\
                    &\le& C\sum_{K\in {\cal T}_h}|S(K)|^{\frac{1}{2}}\|a\|_{W^{2, {\infty}}(K)}\|\nabla {\cal Q}_k^0u_t-\nabla u_t\|_K
                    \|\nabla z_0\|_{K}\nonumber\\
                    &\le& C\sum_{K\in {\cal T}_h}h_{S(K)}\|a\|_{W^{2, {\infty}}(K)}\|\nabla {\cal Q}_k^0u_t-\nabla u_t\|_K
                    \|\nabla z_0\|_{K}\nonumber\\
                    &\le& C\sum_{K\in {\cal T}_h}h_K\|a\|_{W^{2, {\infty}}(K)}\|\nabla {\cal Q}_k^0u_t-\nabla u_t\|_K
                    \|\nabla z_0\|_{K}\nonumber\\
                                      &\le & Ch\|a\|_{2, \infty}Ch^{\lambda_2}\|u_t\|_{\lambda_2+1}\|z_h\|_{1, h}\nonumber\\
                    &\le& Ch^{\lambda_2+1}\|a\|_{2, \infty}\|u_t\|_{\lambda_2+1}\|\rho_t\|,\label{I:1:2}
 		  \end{eqnarray}
   for some non-negative integer $\lambda_2\le k$. In the above estimate, we have used the embedding $H^{2}(K)\hookrightarrow L^{\infty}(K)$ for each $K\in {\cal T}_h$ (cf. Theorem 1.4.6 in \cite{brenner2007mathematical}).
    Set $\lambda=\min\{\lambda_1+1, \lambda_2\}$ and combine above estimates (\ref{k:1})-(\ref{I:1:2}) to obtain
 \begin{eqnarray}
 I_1\le C\|a\|_{2, \infty}h^{\lambda+1}\int_0^T\|u_t\|_{\lambda+1}\|\rho_t\|ds,\;\;0\le \lambda \le k. \label{I:1}
 \end{eqnarray}
Then, following the lines of proof for the Lemma 4.3 in \cite{wang2018systematic}, we obtain
\begin{eqnarray}
 I_2&\le& C\|a\|_{l+1, \infty}h^{\lambda}\int_0^T\|u_t\|_{\lambda+1}\|z-z_h\|_{1, h}ds\nonumber\\
 &\le&C\|a\|_{l+1, \infty}h^{\lambda+1}\int_0^T\|u_t\|_{\lambda+1}\|\rho_t\|ds,\;\;0\le \lambda \le k. \label{I:2}
  \end{eqnarray}
In the last inequality, we have used (\ref{jo:1}).

For $j<l$ and $s=\min\{k, j\}$, Lemma 4.4 in \cite{wang2018systematic} yields
\begin{eqnarray}
l_3(u_t,z_h)&\le& Ch^s\|u_t\|_{s+1}\Bigg(\sum_{K\in {\cal T}_h}\|a\nabla_wz_h\|^2_K\Bigg)^{\frac{1}{2}}\nonumber\\
&\le & Ch^s\|u_t\|_{s+1}\Bigg(\sum_{K\in {\cal T}_h}\|a\|_{L^{\infty}(K)}^2\|\nabla_wz_h\|^2_K\Bigg)^{\frac{1}{2}}\nonumber\\
&\le & Ch^s\|u_t\|_{s+1}Ch\|a\|_{2, \infty}\|z_h\|_{1,h}\nonumber\\
&\le & Ch^{s+1}\|a\|_{2, \infty}\|u_t\|_{s+1}\|\rho_t\|. \label{I:3:1}
\end{eqnarray}
Similarly, for $j\ge l$, we obtain
\begin{eqnarray}
l_3(u_t,z_h)\le Ch^{k+1}\|a\|_{2, \infty}\|u_t\|_{k+1}\|\rho_t\|. \label{I:3:2}
\end{eqnarray}
Combining estimates (\ref{I:3:1})-(\ref{I:3:2}), we have
\begin{equation}\label{I:3}
I_3\le \left\{
  \begin{array}{ll}
    Ch^{s+1}\|a\|_{2, \infty}\int_0^T\|u_t\|_{s+1}\|\rho_t\|ds & \mbox{for}\;j<l, \\~\\
    Ch^{k+1}\|a\|_{2, \infty}\int_0^T\|u_t\|_{k+1}\|\rho_t\|ds, & \mbox{for}\;j\ge l.
  \end{array}
\right.
\end{equation}
Here, $s=\min\{k, j\}$.

Again, for the term $I_4$, we apply Lemma 4.5 and Lemma 4.7 in \cite{wang2018systematic} to have
\begin{eqnarray}
{\cal S}({\cal Q}_hu_t,z_h)&=&{\cal S}({\cal Q}_hu_t,z_h-{\cal Q}_hz)+{\cal S}({\cal Q}_hu_t,{\cal Q}_hz)\nonumber\\
&\le & Ch^s\|u_t\|_{s+1}\|z_h-{\cal Q}_hz\|_{1,h}\nonumber\\&&
+\sKT h^{-1}_K \langle {\cal Q}_m({\cal Q}_k^0u_t-{\cal Q}_j^bu_t), {\cal Q}_m({\cal Q}_k^0z-{\cal Q}_j^bz)\rangle_{\partial K}\nonumber\\
&\le&Ch^{s+1}\|u_t\|_{s+1}\|z\|_2+Ch^s\|u_t\|_2Ch\|z\|_2\nonumber\\
&\le&Ch^{s+1}\|u_t\|_{s+1}\|\rho_t\|, \label{I:4:1}
\end{eqnarray}
with $j<l$ and $s=\{k, j\}.$ Proceeding similarly, for $j\ge l$, we obtain
\begin{eqnarray}
{\cal S}({\cal Q}_hu_t,z_h)
\le Ch^{k+1}\|u_t\|_{k+1}\|\rho_t\|. \label{I:4:2}
\end{eqnarray}
Combining estimates (\ref{I:4:1})-(\ref{I:4:2}), we have
\begin{equation}\label{I:3}
I_4\le \left\{
  \begin{array}{ll}
    Ch^{s+1}\int_0^T\|u_t\|_{s+1}\|\rho_t\|ds & \mbox{for}\;j<l, \\~\\
    Ch^{k+1}\int_0^T\|u_t\|_{k+1}\|\rho_t\|ds & \mbox{for}\;j\ge l,
  \end{array}
\right.
\end{equation}
with $s=\min\{k, j\}$.

Substituting estimates for $I_i\;(1\le i \le 4)$ in (\ref{pa:1}), we obtain
\begin{equation}\label{pa:2}
\|\rho(t)\|^2\le \int_0^T\|\rho_t\|^2ds\le \left\{
  \begin{array}{ll}
    Ch^{2(s+1)}\int_0^T\|u_t\|_{s+1}^2ds & \mbox{for}\;j<l, \\~\\
    Ch^{2(k+1)}\int_0^T\|u_t\|_{k+1}^2 & \mbox{for}\;j\ge l.
  \end{array}
\right.
\end{equation}

Finally, use estimate (\ref{pa:2}) in (\ref{theta:2}) to obtain following $L^{\infty}(L^2)$ norm error estimate.
\begin{them}
 Let $k, j,$ and $l\ge {k-1}$ be the non-negative integers that define the weak finite element space ${\cal V}_h$. Assume that
 $${\cal S}(u_h,v_h)=\sKT h^{-1}_K \langle {\cal Q}_m(u_b-u_0|_{\partial K}), {\cal Q}_m(v_b-v_0|_{\partial K})\rangle_{\partial K},$$
  where ${\cal Q}_m:L^2(\partial K) \to  {\cal P}_m(\partial K)$ is the usual $L^2$ projection operator and $m= \max\{j, l\}$. Then the following error estimates hold ture:
  \begin{description}
  \item {\rm (a)}\;For $j<l$, set $s=\min\{k, j\}$ and assume $s\geq 1$. Assume that the solution of (\ref{1:1})-(\ref{1:2}) is so regular that $u_t\in H^{s+1}(\Omega)\cap H_0^1(\Omega)$. Then
    \begin{eqnarray}
    \|e_h(t)\|^2\le Ch^{2(s+1)}\int_{0}^{t}\|u_t\|_{s+1}^2dt. \label{thm2:a:1}
    \end{eqnarray}
    \item {\rm (b)}\;For $j\ge l$, assume that the solution of (\ref{1:1})-(\ref{1:2}) is so regular that $u_t\in H^{k+1}(\Omega)\cap H_0^1(\Omega)$. Then
    \begin{eqnarray}
    \|e_h(t)\|^2\le Ch^{2(k+1)}\int_{0}^{t}\|u_t\|_{k+1}^2dt. \label{thm2:b:1}
    \end{eqnarray}
  \end{description}
\end{them}

We assume following convergence results for the semidiscrete weak Galerkin approximation with the stabilizer based on element-boundary-discrepancy.
\begin{them}
 Let $k, j,$ and $l\ge {k-1}$ be the non-negative integers that define the weak finite element space ${\cal V}_h$. Assume that
 $${\cal S}(u_h,v_h)=
    \sKT h^{-1}_K \langle u_b-u_0|_{\partial K}, v_b-v_0|_{\partial K}\rangle_{\partial K}.$$
   Then, we have following error estimates
    \begin{eqnarray}
    \|e_h(t)\|+h\|e_h(t)\|_{1, h}\le Ch^{s+1}\Bigg(\int_{0}^{T}\|u_t\|_{s+1}^2dt\Bigg)^{\frac{1}{2}}, \label{thm3}
    \end{eqnarray}
    where $s=\min\{k, j\}$.
\end{them}

  \section{Discrete time WG Finite Element Method}\label{s6}\se
   We now turn our attention to some discrete time weak Galerkin procedures.
A discrete-in-time scheme based on backward Euler method for approximating exact solution $u$ is discussed in this section. Optimal pointwise-in-time error estimate in both discrete $H^1$ and $L^2$ norms are established.

First we divide the time interval $J=[0, T]$ into $M$ equally spaced subintervals $I_n=(t_{n-1}, t_n]$, $n=1,2,\ldots,M$ with $t_0=0$, and $t_M=T$ and $\tau=t_n-t_{n-1}$, the time step.
For a sequence $\{p^n\}_{n=0}^M \subset L^2(\Omega)$, we define
\begin{eqnarray*}
\partial_\tau p^n=\frac{p^{n}-p^{n-1}}{\tau},\;\;\;n=1,\ldots,M.
\end{eqnarray*}
Also, for a continuous mapping $\phi:[0, T]\rightarrow L^2(\Omega)$, we define $\phi^n=\phi(., t_n)$, $ 0\leq n\leq M$. 

With the above notation, we now introduce the fully discrete weak Galerkin finite element approximation to the problem (\ref{1:1})-(\ref{1:2}): Let $U_h^0={\cal R}_h\psi$ and $U_h^n=\{U_0^n, U_b^n\}\in V_h^0$
be the fully discrete solution of $u$ at $t=t_n$ which we shall define through the following scheme
\begin{equation}
(\partial_\tau U_h^n, v_h) + {\cal A}(U_h^n, v_h) = ({f}^{n}, v_h)\;\;\forall
v_h\in V_{h}^0,\;\;n=1,\ldots,M. \label{4:1}
\end{equation}

For fully discrete error estimates, we now split the errors at $t=t_n$ as follows
\begin{eqnarray*}
{u}^n-{U}_h^n={u}^n-{\cal Q}_h{u}^n + {\cal Q}_h{u}^n - {U}_h^n.
\end{eqnarray*}
We denote our error as
\begin{eqnarray*}
{e}^n = {U}_h^n-{\cal Q}_h{u}^n=\{{e}_0^n,{ e}_b^n\}.
\end{eqnarray*}

\subsection{\normalsize Error estimates with  projected element-boundary-discrepancy}
 Convergence results for the fully discrete weak Galerkin approximation with the stabilizer based on projected element-boundary-discrepancy
 are presented.

Using $\rho$ and $\theta$, error $e^n$ can be further separated as
\begin{eqnarray}
{e}^n =\theta^n+\rho^n, \label{4:2}
\end{eqnarray}
where $\theta^n=U_h^n-{\cal R}_hu^n$ and $\rho^n={\cal R}_hu^n-{\cal Q}_hu^n$.

For $\theta^n$, we have the following error equation
\begin{eqnarray}
({\partial}_\tau\theta^n,v_0)+{\cal A}(\theta^n,v)&=&-({\partial_\tau}{\cal R}_hu^n-u_t^n,v_0)\nonumber\\&:=&-(w^n,v_0)\ \forall {v}=\{v_0,v_b\}\in {V}_h^0,\label{4:11}
\end{eqnarray}
where $w^n={\partial_\tau}{\cal R}_hu^n-u_t^n$.
For simplicity of the exposition,
we write $w^n={\cal R}_1^n+{\cal R}_2^n$, where
${\cal R}_1^n={\partial_\tau}{\cal R}_hu^n-{\partial_\tau}u^n$ and ${\cal R}_2^n={\partial_\tau}u^n-u_t^n$.

Set $v=\theta^n$ in (\ref{4:11}), we have
\begin{eqnarray}
({\partial}_\tau\theta^n,\theta^n) +{\cal A}(\theta^n,\theta^n) \leq \|w^n\|\|\theta^n\|. \nonumber
\end{eqnarray}
Then using the positivity of ${\cal A}(\cdot, \cdot)$, we obtain
$$\|\theta^n\|\leq \|\theta^{n-1}\|+ \tau\|w^n\|,$$
where we have used the fact that $\theta^0=U_h^0-{\cal R}_hu^0=0.$ Hence, we have
\begin{eqnarray}
\|\theta^n\|\leq \tau\sum_{j=1}^{n}\|w^j\|\leq \tau\sum_{j=1}^{n}\|{\cal R}_1^j\|+\tau\sum_{j=1}^{n}\|{\cal R}_2^j\|. \label{4:12}
\end{eqnarray}
For the term ${\cal R}_1^j$, it is easy to verify that
\begin{eqnarray*}
\tau{\cal R}_1^j=\int_{t_{j-1}}^{t_j}({\cal R}_hu_t-u_t)ds,\label{n:c:1}
\end{eqnarray*}
which together with estimates (\ref{pa:2}) leads to the following
\begin{eqnarray}
 \tau\sum_{j=1}^{n}\|{\cal R}_1^j\|  \leq
  \left\{
  \begin{array}{ll}
  Ch^{(s+1)}\Big(\int_{0}^{T}\|u_t\|_{s+1}^2ds\Big)^{\frac{1}{2}} & \mbox{for}\;j<l, \\~\\
  Ch^{k+1}\Big(\int_{0}^{T}\|u_t\|_{k+1}^2ds\Big)^{\frac{1}{2}} & \mbox{for}\;j\ge l.
  \end{array}
  \right.   \label{4:14}
\end{eqnarray}
Now, for the term ${\cal R}_2$, we use Taylor's series expansion to have
\begin{eqnarray}
 \tau\sum_{j=1}^{n}\|{\cal R}_2^j\|\leq C{\tau}\int_{0}^{T}\|u_{tt}\|ds. \label{4:13}
\end{eqnarray}

Finally, estimates (\ref{4:12})-(\ref{4:13}) together with (\ref{pa:2}) leads to following $L^2$ norm error estimates.
\begin{them}  \label{fully:thm:1}
	Let $k,\; j$ and $l\ge {k-1}$ be the non-negative integers that define the weak finite element space ${\cal V}_h$. Assume that
	$${\cal S}(u_h,v_h)=\sKT h^{-1}_K \langle {\cal Q}_m(u_b-u_0|_{\partial K}), {\cal Q}_m(v_b-v_0|_{\partial K})\rangle_{\partial K},$$
	where ${\cal Q}_m:L^2(\partial K) \to  {\cal P}_m(\partial K)$ is the usual $L^2$ projection operator and $m=\max\{j, l\}$. Then the following error estimates hold true:
	\begin{description}
		\item {\rm (a)}\;For $j<l$, set $s=\min\{k, j\}$ and assume $s\geq 1$. Assume that the solution of (\ref{1:1})-(\ref{1:2}) is so regular that $u_t\in H^{s+1}(\Omega)\cap H_0^1(\Omega)$. Then
		\begin{eqnarray}
		\|e^n\|^2\le C\Big(h^{2(s+1)}+\tau^2\Big)\int_{0}^{T}\Big(\|u_t\|_{s+1}^2+\|u_{tt}\|^2\Big)dt. \label{thm4:a:1}
		\end{eqnarray}
		\item {\rm (b)}\;For $j\ge l$, assume that the solution of (\ref{1:1})-(\ref{1:2}) is so regular that $u_t\in H^{k+1}(\Omega)\cap H_0^1(\Omega)$. Then
		\begin{eqnarray}
		\|e^n\|^2\le C\Big(h^{2(k+1)}+\tau^2\Big)\int_{0}^{T}\Big(\|u_t\|_{k+1}^2+\|u_{tt}\|^2\Big)dt. \label{thm4:b:1}
		\end{eqnarray}
	\end{description}
\end{them}

Next, setting $v_h=-\tau\partial_\tau \theta^n$ in $(\ref{4:11})$, we have
 \begin{equation}
\tau\|\partial_\tau \theta^n\|^2+{\cal A}(\theta^{n},\theta^n-\theta^{n-1})\leq
\tau\|w^n\|\|\partial_\tau\theta^n\|,\label{:1:tc}
\end{equation}
which yields
\begin{eqnarray}
\tau\|\partial_\tau \theta^n\|^2+\tnorm{\theta^n}^2-\tnorm{\theta^{n-1}}^2&\leq&
C\tau\|w^n\|^2\nonumber\\
&\le& C\tau(\|{\cal R}_1^n\|^2+\|{\cal R}_2^n\|^2).\label{c:1:tc}
\end{eqnarray}
From (\ref{n:c:1}), we note that
\begin{eqnarray*}
\tau{\cal R}_1^j\le \tau^{\frac{1}{2}}\Bigg(\int_{t_{j-1}}^{t_j}({\cal R}_hu_t-u_t)^2dt\Bigg)^{\frac{1}{2}}
\end{eqnarray*}
so that following estimates hold true
\begin{eqnarray}
 \tau\sum_{j=1}^{n}\|{\cal R}_1^j\|^2  \leq
  \left\{
  \begin{array}{ll}
  Ch^{2(s+1)}\int_{0}^{T}\|u_t\|_{s+1}^2ds & \mbox{for}\;j<l, \\~\\
  Ch^{2(k+1)}\int_{0}^{T}\|u_t\|_{k+1}^2ds & \mbox{for}\;j\ge l.
  \end{array}
  \right.   \label{n:4:14}
\end{eqnarray}
Again, we know that
\begin{eqnarray*}
{\cal R}_2^j =\frac{u^j-u^{j-1}}{{\tau}}-u_t^{j}=-\frac{1}{{\tau}}\int_{t_{j-1}}^{t_j}(s-t_{j-1})u_{tt}ds.
\end{eqnarray*}
Hence, we have
\begin{eqnarray*}
|{\cal R}_2^j|^2&\le &\frac{1}{{\tau^2}}\Bigg(\int_{t_{j-1}}^{t_j}(s-t_{j-1})^2ds\Bigg)\Bigg(\int_{t_{j-1}}^{t_j}u_{tt}^2ds\Bigg)\\
&\le& C\tau\int_{t_{j-1}}^{t_j}u_{tt}^2ds,
\end{eqnarray*}
which integration over $\Omega$ yields
\begin{eqnarray}
\|{\cal R}_1^j\|^2
\leq C{\tau}\int_{t_{j-1}}^{t_j}\|u_{tt}\|^2ds.\label{n:4:6}
\end{eqnarray}

Now, estimates (\ref{c:1:tc})-(\ref{n:4:6}) together with (\ref{pa:2}) leads to following discrete $H^1$ norm error estimates.
\begin{them} \label{fully:thm:2}
	Let $k,\; j$ and $l\ge {k-1}$ be the non-negative integers that define the weak finite element space ${\cal V}_h$. Assume that
	$${\cal S}(u_h,v_h)=\sKT h^{-1}_K \langle {\cal Q}_m(u_b-u_0|_{\partial K}), {\cal Q}_m(v_b-v_0|_{\partial K})\rangle_{\partial K},$$
	where ${\cal Q}_m:L^2(\partial K) \to  {\cal P}_m(\partial K)$ is the usual $L^2$ projection operator and $m=\max\{j, l\}$. Then the following error estimates hold true:
	\begin{description}
		\item {\rm (a)}\;For $j<l$, set $s=\min\{k, j\}$ and assume $s\geq 1$. Assume that the solution of (\ref{1:1})-(\ref{1:2}) is so regular that $u_t\in H^{s+1}(\Omega)\cap H_0^1(\Omega)$. Then
		\begin{eqnarray}
		\tnorm{e^n}^2\le C\Big(h^{2s}+\tau^2\Big)\int_{0}^{T}\Big(\|u_t\|_{s+1}^2+\|u_{tt}\|^2\Big)dt. \label{thm3:a:1}
		\end{eqnarray}
		\item {\rm (b)}\;For $j\ge l$, assume that the solution of (\ref{1:1})-(\ref{1:2}) is so regular that $u_t\in H^{k+1}(\Omega)\cap H_0^1(\Omega)$. Then
		\begin{eqnarray}
		\tnorm{e^n}^2\le C\Big(h^{2k}+\tau^2\Big)\int_{0}^{T}\Big(\|u_t\|_{k+1}^2+\|u_{tt}\|^2\Big)dt. \label{thm3:b:1}
		\end{eqnarray}
	\end{description}
\end{them}

We assume following convergence results for the fully discrete weak Galerkin approximation with the stabilizer based on element-boundary-discrepancy.
\begin{them}  \label{fully:thm:3}
	Let $k, j,$ and $l\ge {k-1}$ be the non-negative integers that define the weak finite element space ${\cal V}_h$.  Assume that
	$${\cal S}(u_h,v_h)=
    \sKT h^{-1}_K \langle u_b-u_0|_{\partial K}, v_b-v_0|_{\partial K}\rangle_{\partial K}.$$
	Then, we have following error estimates
		\begin{eqnarray}
		\|e^n\|^2&\le& C\Big(h^{2(s+1)}+\tau^2\Big)\int_{0}^{T}\Big(\|u_t\|_{s+1}^2+\|u_{tt}\|^2\Big)dt, \label{thm4:bb:1}\\
			\tnorm{e^n}^2&\le& C\Big(h^{2s}+\tau^2\Big)\int_{0}^{T}\Big(\|u_t\|_{s+1}^2+\|u_{tt}\|^2\Big)dt, \label{thm4:bb:1}
		\end{eqnarray}
	where $s=\min\{k, j\}$.
	\end{them}

\begin{remm}
 From Theorems \ref{fully:thm:1}-\ref{fully:thm:3} it is clear that the method of projected element-boundary-discrepancy is more accurate than the method of element-boundary-discrepancy
 for the case of $j=k-1$ and $l=k-1.$  For the numerical validation, we refer to Tables \ref{table1}-\ref{table8}.
\end{remm}

 \section{Numerical Experiments}\label{sec7}\se
In this section we will explore the results of computations for the parabolic problems (\ref{1:1})-(\ref{1:2}) in $\Omega\times J$, where $\Omega=(0, 1)\times (0, 1)$ and $J=[0,1]$
with selected values on the degree of polynomials in the weak Galerkin finite element space.
The coefficient matrix $a$ is given by identity matrix $I$, the load function $f$, initial data $\psi$ and the Dirichlet boundary value are selected in such a way that exact solution is
$u=\exp(-t)\sin(\pi x)\sin(\pi y)$. In this test problem, triangular mesh is used. We have done uniform partitioning of the domain into $n \times n$ sub rectangles which is followed by dividing each rectangular element by the diagonal line with mesh size $h =1/n,$ where $n$ is any non-negative integer. 
Further, we set $\tau=\mbox{O}(h^{\gamma+1})$, where $\gamma$ is selected according to Theorems \ref{fully:thm:1}-\ref{fully:thm:3} so that optimal order of convergence is maintained.

Let $U_h^n$ be the weak Galerkin solution defined by (\ref{4:1}). Then, we have calculated the following error
\begin{eqnarray*}
	{e}^n = {U}_h^n-{\cal Q}_h{u}^n=\{{ e}_0^n, {e}_b^n\},
\end{eqnarray*}
with respect to triple bar norm and the $L^2$ norm at final time $T=1$.

\begin{center}
       	\begin{figure}[t]
       		\includegraphics[width=5.5cm, height=4.5cm]{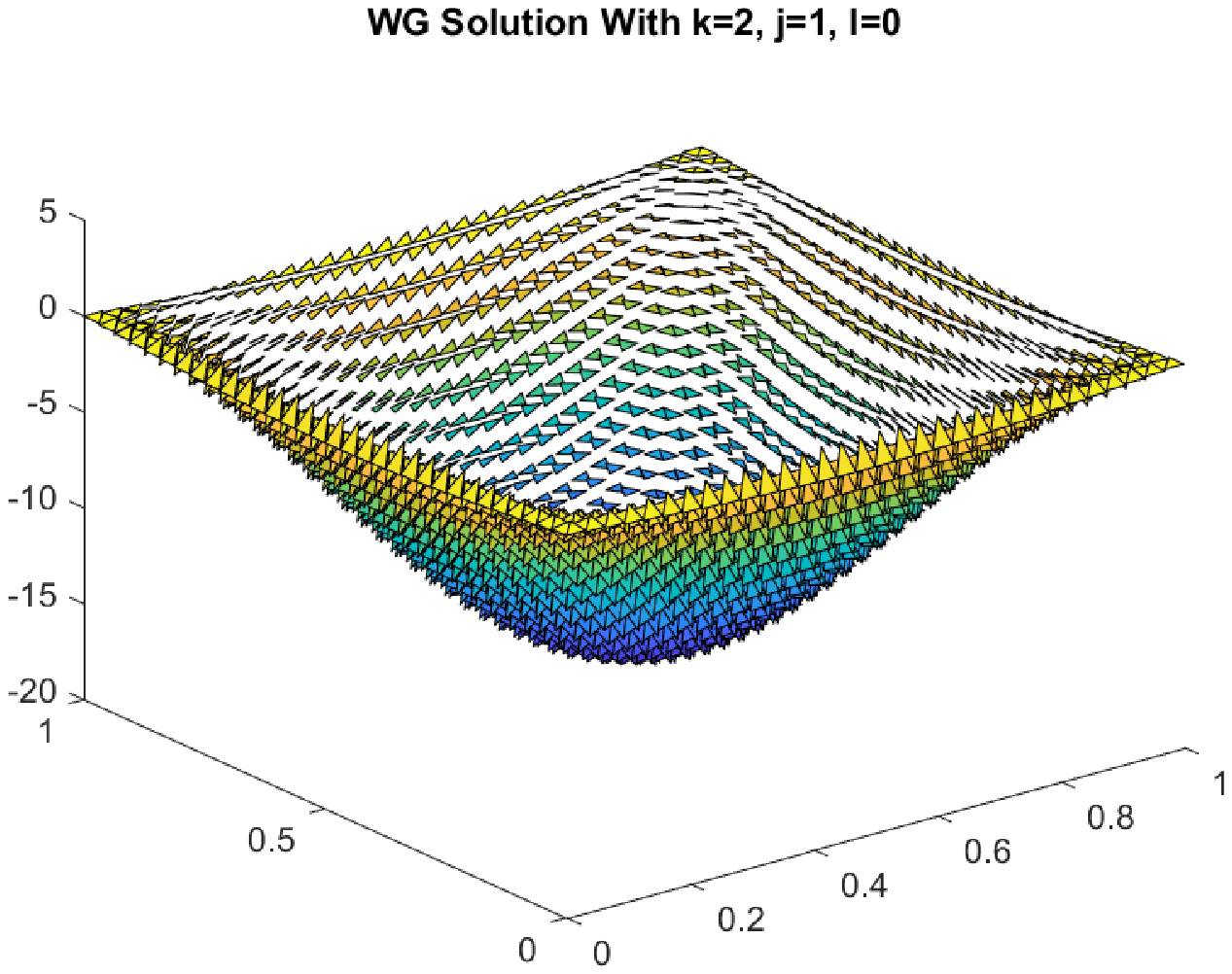}\;\;\;\includegraphics[width=5.5cm, height=4.5cm]{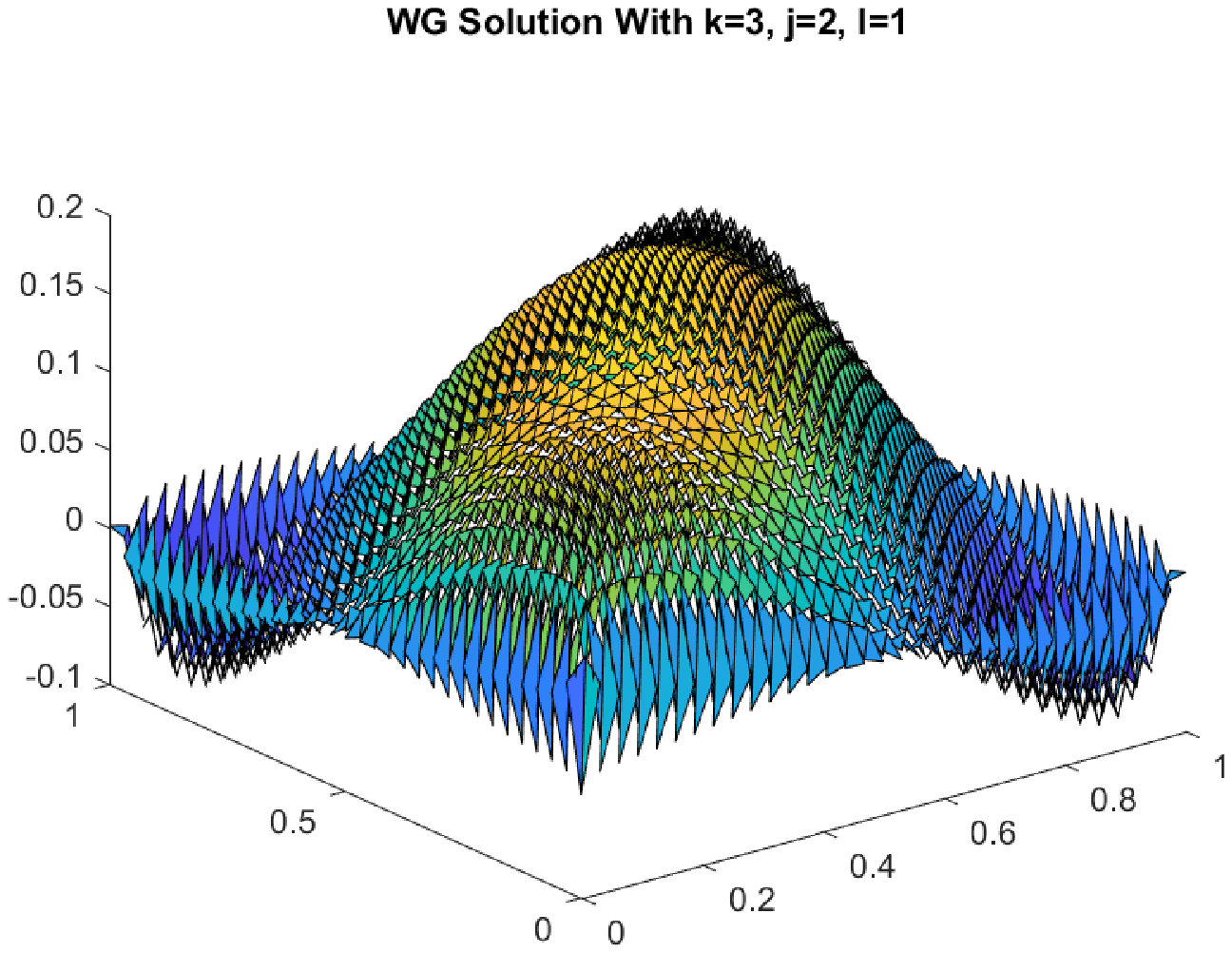}\;\;
       		\includegraphics[width=5.5cm,height=4.5cm]{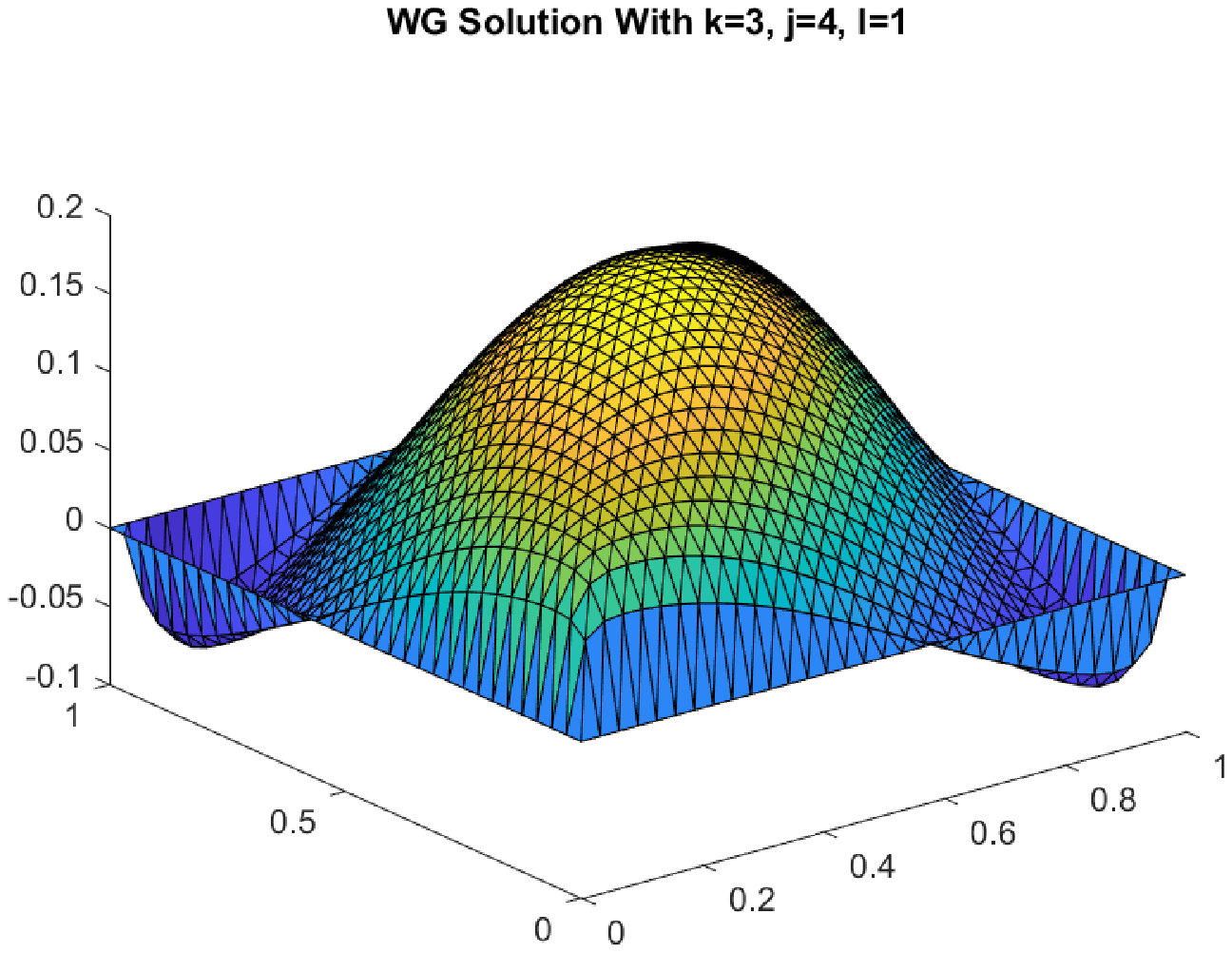}\;\;\;\includegraphics[width=5.5cm,height=4.5cm]{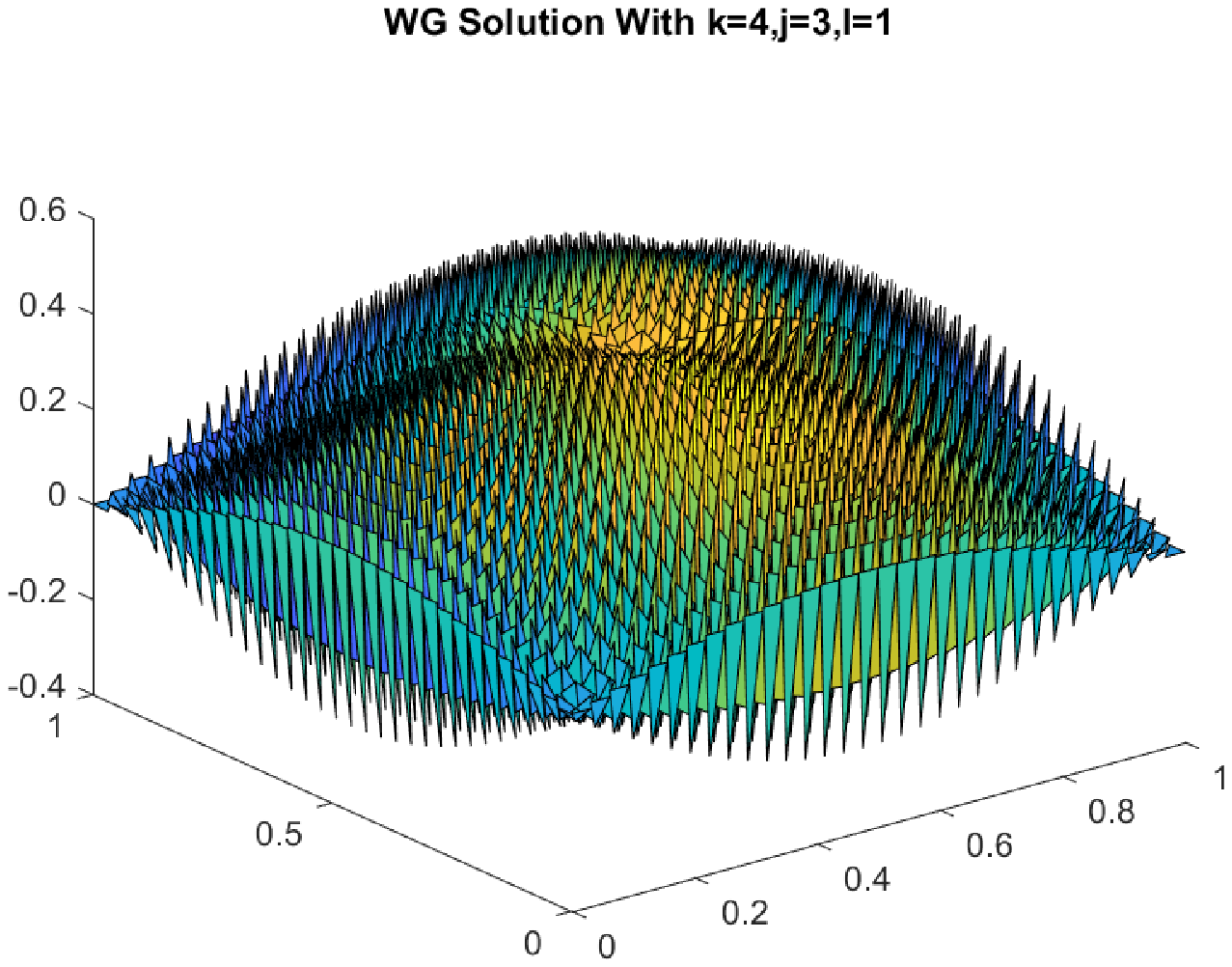}\centering
       		\caption{Plots of the WG approximations at time $t=1$ for the method of projected element-boundary-discrepancy with $h=1/32.$}\label{fig5.1}
       	\end{figure}
       \end{center}

Recall that the stabilizer for the method of projected element-boundary-discrepancy is given by  $${\cal S}(u_h,v_h)=\sKT h^{-1}_K \langle {\cal Q}_m(u_b-u_0|_{\partial K}), {\cal Q}_m(v_b-v_0|_{\partial K})\rangle_{\partial K},$$
where $m=\min\{j, l\}.$ For different values of $k\;(1\le k\le 4),\;j\;(0\le j \le 4)$ and $l\;(0\le l \le 4),$
we have implemented the corresponding WG scheme (\ref{4:1}) for the problem (\ref{1:1})-(\ref{1:2}).
The rate of convergence for each combination is reported in Tables \ref{table1}-\ref{table4}, where $\mbox{NI}$
means the corresponding WG scheme is unstable or not consistent. Detailed computational data can be found in the Appendix. The convergence order for each particular combination is indicated in the form $n/m$,
where $n$ stand for the order of convergence in the triple bar norm and $m$ for the order of convergence in the $L^2$ norm. For example, $2/3$ would mean that the method is convergent at the rate of $h^2$ in the triple bar
norm and $h^3$ in the $L^2$ norm. For $l<k-1$, the method works poorly. This is an observation from the computation. For instance, we refer to Figure \ref{fig5.1}.
	
\begin{table}[h]
	\caption{ $~~~~$ Order of convergence for $k$=1 with stabilizer term\newline
		 $~~~~~~~~~~~~~~~~~~~~ h^{-1}_K \langle {\cal Q}_m(u_b-u_0|_{\partial K}),{\cal Q}_m(v_b-v_0|_{\partial K})\rangle$ } \centering
	 \vspace{0.3cm}
	\begin{tabular}{|c|c|c|c|c|c|}
		\hline
		$k=1$ & $ j=0 $& $ j=1 $& $j=2$  &$ j=3 $& $j=4$ \\ \hline
		$l=0$  & 1/2 & 1/2 &  1/2  & 1/2 & 1/2  \\
		$l=1$  & 0/0 & 1/2 &  1/2  & 1/2 & 1/2     \\
		$l=2$  & 0/0 & 1/2 &  1/2  & 1/2 & 1/2  \\
		$l=3$  & 0/0 & 1/2 &  1/2  & 1/2 & 1/2    \\
		$l=4$  & 0/0 & 1/2 &  1/2  & 1/2 & 1/2     \\ \hline
	\end{tabular}
\label{table1}
		\caption{ $~~~~$ Order of convergence for $k$=2 with stabilizer term\newline
		$~~~~~~~~~~~~~~~~~~~~ h^{-1}_K \langle {\cal Q}_m(u_b-u_0|_{\partial K}),{\cal Q}_m(v_b-v_0|_{\partial K})\rangle$ } \centering
	\vspace{0.3cm}
	\begin{tabular}{|c|c|c|c|c|c|}
		\hline
		$k=2$ & $ j=0 $& $ j=1 $& $j=2$  &$ j=3 $& $j=4$ \\ \hline
		$l=0$  & \mbox{NI}  & \mbox{NI} &  \mbox{NI}  & \mbox{NI} & \mbox{NI}  \\
		$l=1$  & 0/0  & 2/3 &  2/3  & 2/3 & 2/3     \\
		$l=2$  & 0/0  & 1/2 &  2/3  & 2/3 & 2/3  \\
		$l=3$  & 0/0  & 1/2 &  2/3  & 2/3 & 2/3  \\
		$l=4$  & 0/0  & 1/2 &  2/3  & 2/3 & 2/3  \\  \hline
	\end{tabular}
	\label{table2}
		\caption{ $~~~~$ Order of convergence for $k$=3 with stabilizer term\newline
		$~~~~~~~~~~~~~~~~~~~~ h^{-1}_K \langle {\cal Q}_m(u_b-u_0|_{\partial K}),{\cal Q}_m(v_b-v_0|_{\partial K})\rangle$ } \centering
	\vspace{0.3cm}
	\begin{tabular}{|c|c|c|c|c|c|}
		\hline
		$k=3$ & $ j=0 $& $ j=1 $& $j=2$  &$ j=3 $& $j=4$ \\ \hline
		$l=0$ &  \mbox{NI}  & \mbox{NI} &  \mbox{NI}  & \mbox{NI}& \mbox{NI}  \\
		$l=1$ &  \mbox{NI}  & \mbox{NI} &  \mbox{NI}  & \mbox{NI} &\mbox{NI}  \\
		$l=2$ &  0/0  & 1/2 &  3/4  & 3/4 & 3/4  \\
		$l=3$ &  0/0  & 1/2 &  2/3  & 3/4 & 3/4  \\
		$l=4$ &  0/0  & 1/2 &  2/3  & 3/4 & 3/4  \\ \hline
	\end{tabular}
	\label{table3}
		\caption{ $~~~~$ Order of convergence for $k$=4 with stabilizer term\newline
		$~~~~~~~~~~~~~~~~~~~~ h^{-1}_K \langle {\cal Q}_m(u_b-u_0|_{\partial K}),{\cal Q}_m(v_b-v_0|_{\partial K})\rangle$ } \centering
	\vspace{0.3cm}
	\begin{tabular}{|c|c|c|c|c|c|}
		\hline
		$k=4$ & $ j=0 $& $ j=1 $& $j=2$  &$ j=3 $& $j=4$ \\ \hline
		$l=0$ &  \mbox{NI}  & \mbox{NI} &  \mbox{NI}  & \mbox{NI} & \mbox{NI}  \\
		$l=1$ &  \mbox{NI}  & \mbox{NI} &  \mbox{NI}  & \mbox{NI} & \mbox{NI}  \\
		$l=2$ &  \mbox{NI}  & \mbox{NI} &  \mbox{NI}  & \mbox{NI} & \mbox{NI}  \\
		$l=3$ &  0/0  & 1/2 &  2/3  & 4/5 & 4/5  \\
		$l=4$ & 0/0  & 1/2 &  2/3  & 3/4 & 4/5  \\ \hline
	\end{tabular}
	\label{table4}
\end{table}

\begin{center}
\begin{figure}[h]
\includegraphics[width=5.5cm, height=4.5cm]{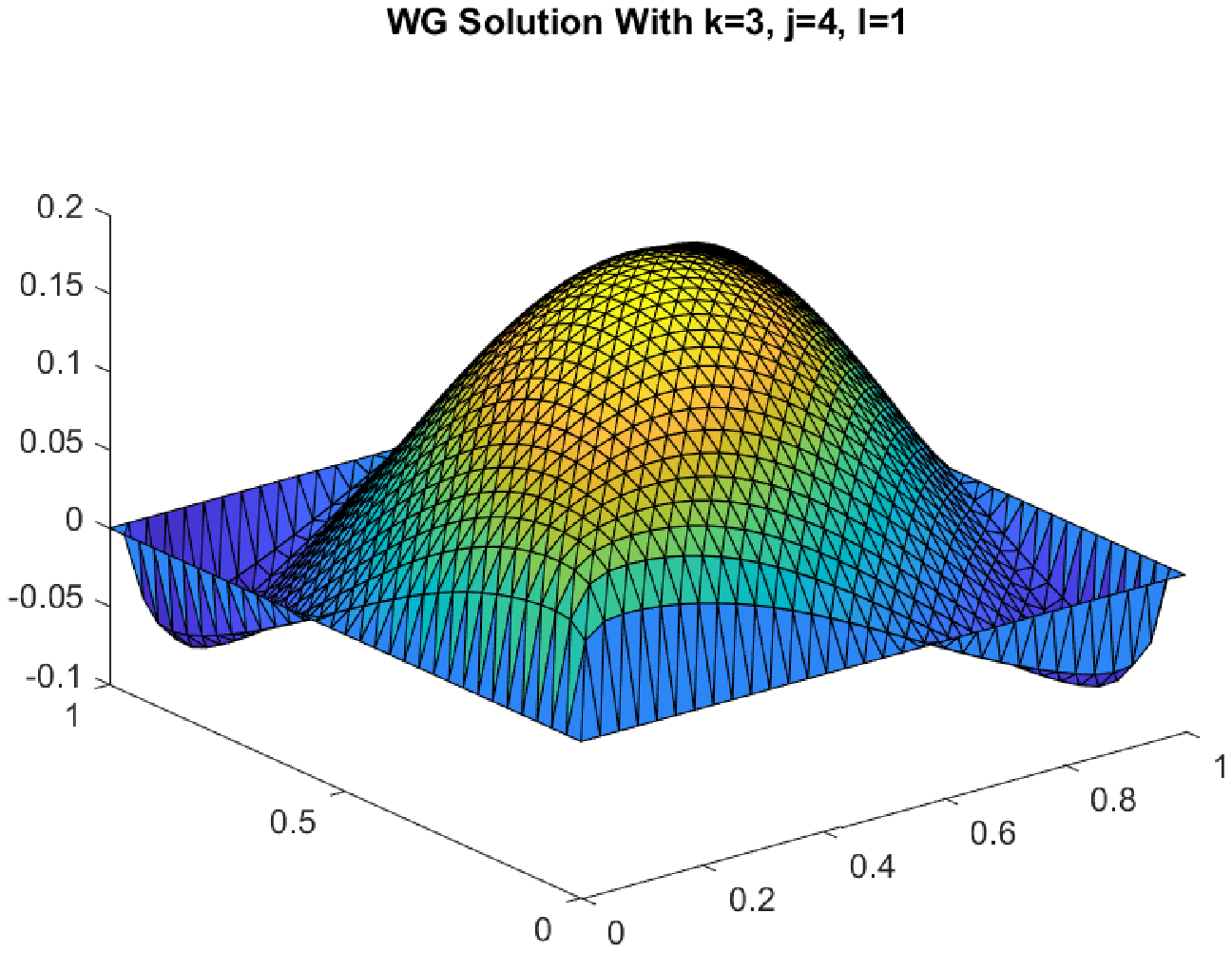}\;\;\;\includegraphics[width=5.5cm, height=4.5cm]{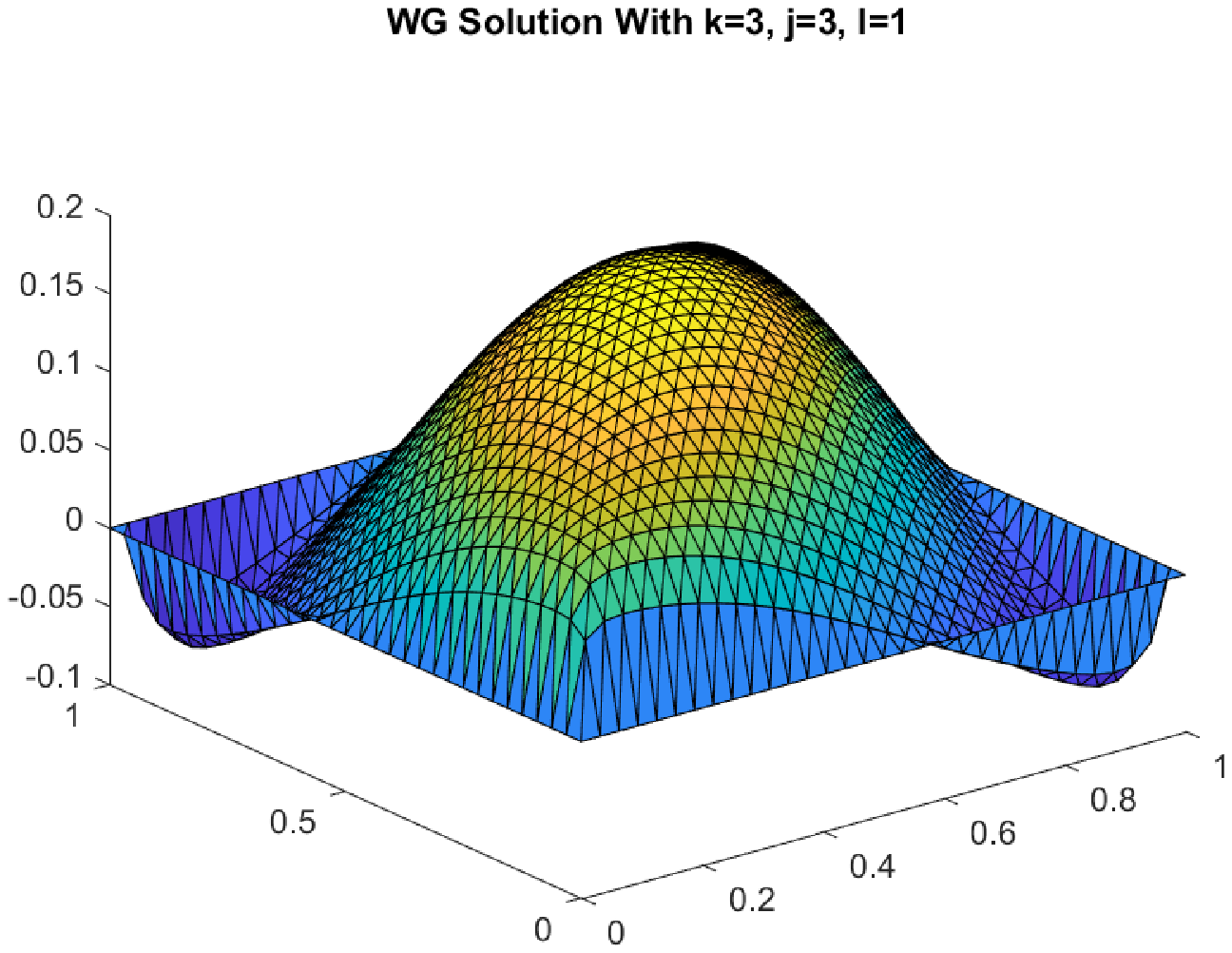}\;\;
\includegraphics[width=5.5cm,height=4.5cm]{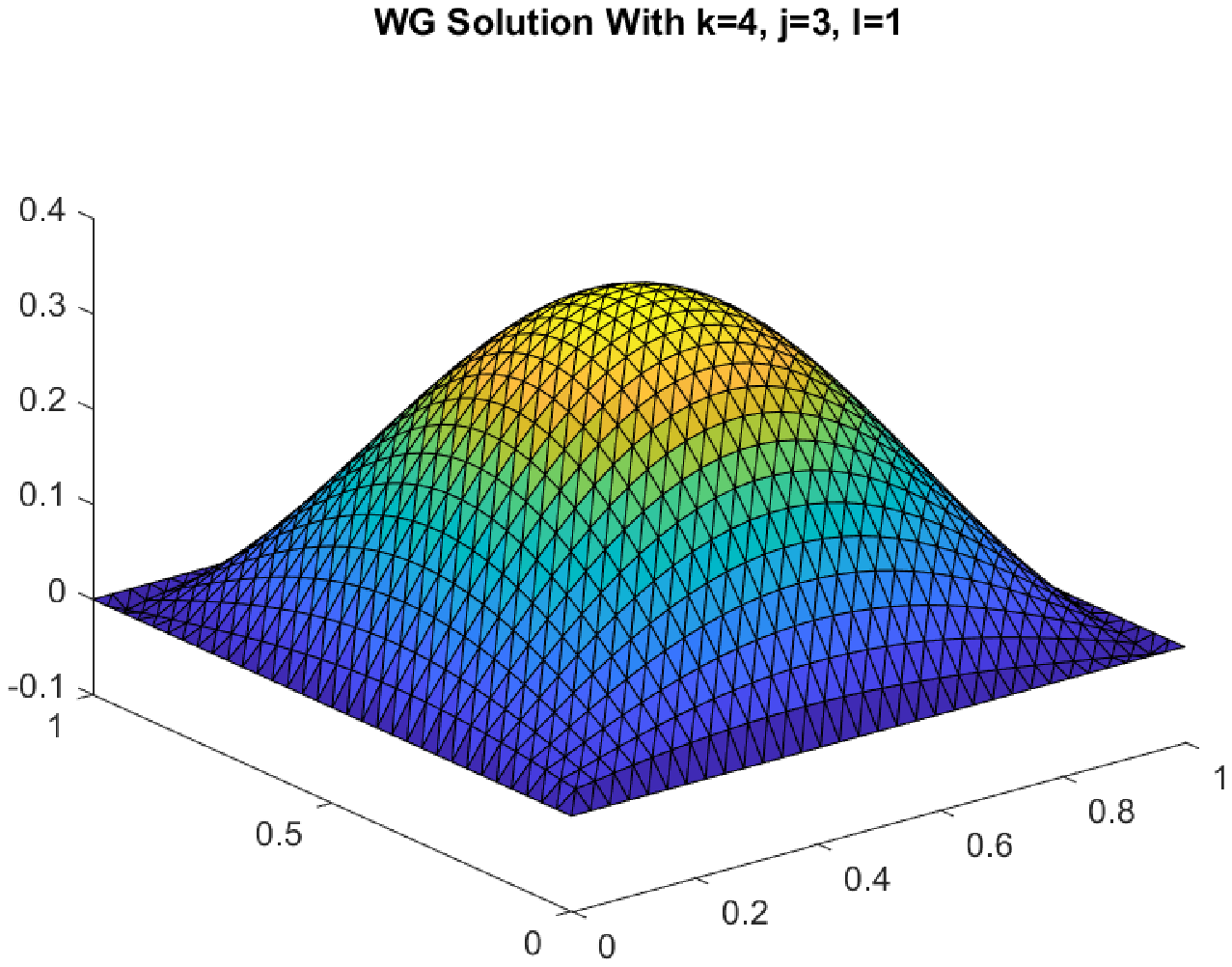}\;\;\;\includegraphics[width=5.5cm,height=4.5cm]{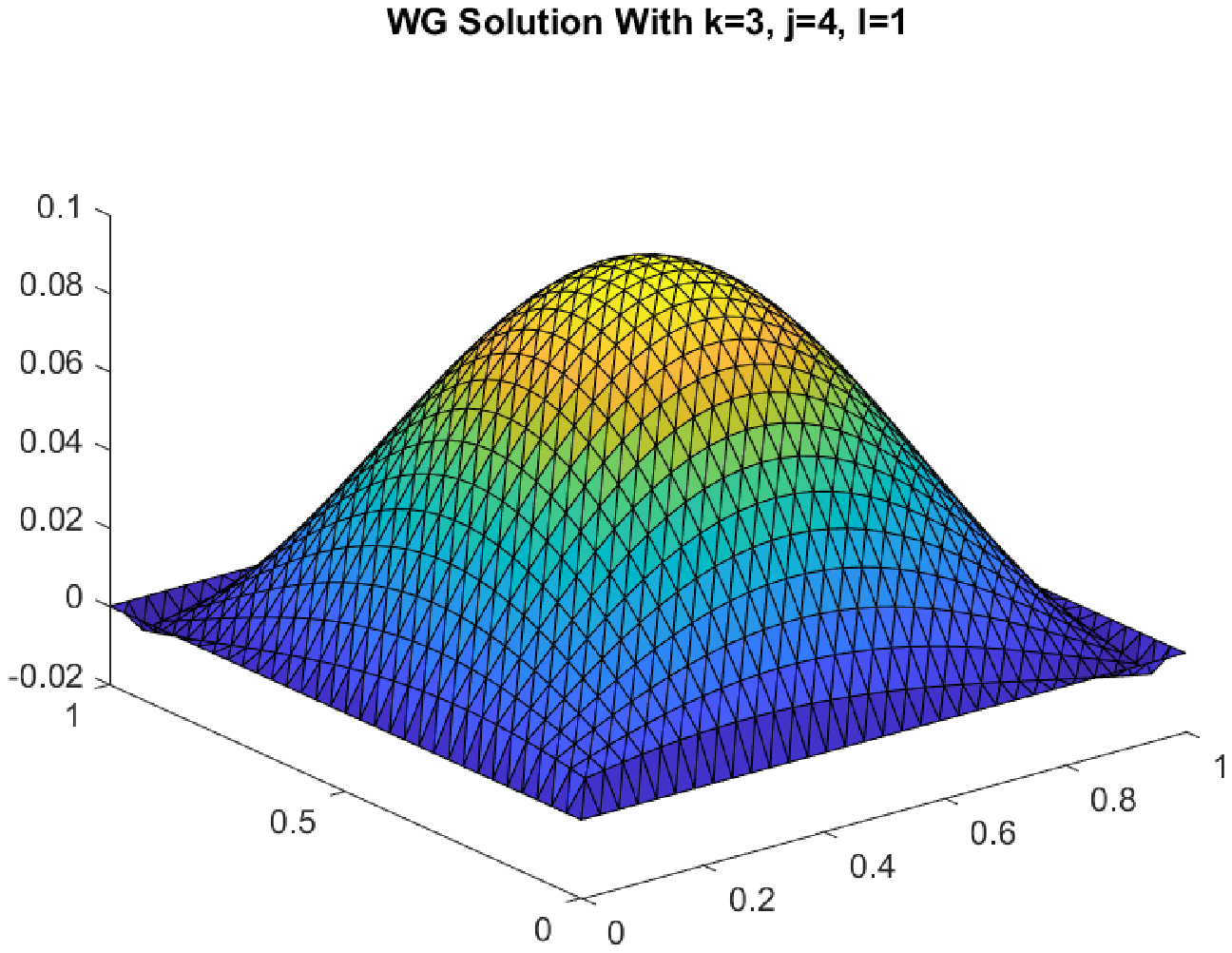}\centering
\caption{Plots of the WG approximations at time $t=1$ for the method of element-boundary-discrepancy with $h=1/32.$}\label{fig5.2}
\end{figure}
\end{center}

The method of element-boundary-discrepancy is based on the selection of stabilizer ${\cal S}(u_h,v_h)=\sKT h^{-1}_K \big\langle u_b-u_0|_{\partial K}, v_b-v_0|_{\partial K} \big\rangle_{\partial K}.$
For all the values of $k =1,\ldots,4, j=0,\ldots,4$, and $l=0,\ldots,4$, we have implemented the corresponding WG finite element scheme (\ref{4:1}). The order of convergence for each combination is listed in
Tables \ref{table5}-\ref{table8}.
Tables \ref{table5}-\ref{table8} suggest that the WG algorithms corresponding to the case of $l=k-2$ and $j<k$ are solvable. For $l=k-2$ with $j\ge k$ and $l<k-2$, the method is solvable but not consistent.
At present, we do not have mathematical justification. This is an observation from the numerical experiments, which is illustrated in Figure \ref{fig5.2}.

\begin{table}[h]
		\caption{ $~~~~$ Order of convergence for $k$=1 with stabilizer term\newline
			$~~~~~~~~~~~~~~~~~~~~ h^{-1}_K \langle u_b-u_0|_{\partial K},v_b-v_0|_{\partial K}\rangle$ } \centering
	\vspace{0.3cm}
	\begin{tabular}{|c|c|c|c|c|c|}
		\hline
		$k=1$ & $ j=0 $& $ j=1 $& $j=2$  &$ j=3 $& $j=4$ \\ \hline
		$l=0$  & 0/0 & 1/2 &  1/2  & 1/2 & 1/2  \\
		$l=1$  & 0/0 & 1/2 &  1/2  & 1/2 & 1/2     \\
		$l=2$  & 0/0 & 1/2 &  1/2  & 1/2 & 1/2  \\
		$l=3$  & 0/0 & 1/2 &  1/2  & 1/2 & 1/2    \\
		$l=4$  & 0/0 & 1/2 &  1/2  & 1/2 & 1/2     \\ \hline
	\end{tabular}
	\label{table5}
\end{table}

\begin{table}[h]
\caption{ $~~~~$ Order of convergence for $k$=2 with stabilizer term\newline
	$~~~~~~~~~~~~~~~~~~~~ h^{-1}_K \langle u_b-u_0|_{\partial K},v_b-v_0|_{\partial K}\rangle$ } \centering
\vspace{0.3cm}
	\begin{tabular}{|c|c|c|c|c|c|}
		\hline
		$k=2$ & $ j=0 $& $ j=1 $& $j=2$  &$ j=3 $& $j=4$ \\ \hline
		$l=0$    & 0/0  & 1/2 &  \mbox{NI}  & \mbox{NI} & \mbox{NI}  \\
		$l=1$  & 0/0  & 1/2 &  2/3  & 2/3 & 2/3     \\
		$l=2$   & 0/0  & 1/2 &  2/3  & 2/3 & 2/3  \\
		$l=3$   & 0/0  & 1/2 &  2/3  & 2/3 & 2/3  \\
		$l=4$  & 0/0  & 1/2 &  2/3  & 2/3 & 2/3  \\  \hline
	\end{tabular}
	\label{table6}
\end{table}

\begin{table}[h]
\caption{ $~~~~$ Order of convergence for $k$=3 with stabilizer term\newline
	$~~~~~~~~~~~~~~~~~~~~ h^{-1}_K \langle u_b-u_0|_{\partial K},v_b-v_0|_{\partial K}\rangle$ } \centering
\vspace{0.3cm}
	\begin{tabular}{|c|c|c|c|c|c|}
		\hline
		$k=3$ & $ j=0 $& $ j=1 $& $j=2$  &$ j=3 $& $j=4$ \\ \hline
		$l=0$ &  \mbox{NI}   & \mbox{NI}  & \mbox{NI}   & \mbox{NI} & \mbox{NI}   \\
		$l=1$ &  0/0  & 1/2 &  2/3  & \mbox{NI} & \mbox{NI}  \\
		$l=2$ &  0/0  & 1/2 &  2/3  & 3/4 & 3/4  \\
		$l=3$ &  0/0  & 1/2 &  2/3  & 3/4 & 3/4  \\
		$l=4$ & 0/0  & 1/2 &  2/3  & 3/4 & 3/4  \\ \hline
	\end{tabular}
	\label{table7}
	\caption{ $~~~~$ Order of convergence for $k$=4 with stabilizer term\newline
	$~~~~~~~~~~~~~~~~~~~~ h^{-1}_K \langle u_b-u_0|_{\partial K},v_b-v_0|_{\partial K}\rangle$ } \centering
\vspace{0.3cm}
	\begin{tabular}{|c|c|c|c|c|c|}
		\hline
		$k=4$ & $ j=0 $& $ j=1 $& $j=2$  &$ j=3 $& $j=4$ \\ \hline
		$l=0$ &  \mbox{NI}  & \mbox{NI} &  \mbox{NI}  & \mbox{NI} & \mbox{NI}  \\
		$l=1$ &  \mbox{NI}  & \mbox{NI} &  \mbox{NI}  & \mbox{NI} & \mbox{NI}  \\
		$l=2$ &  0/0  & 1/2 &  2/3  & 3/4 & \mbox{NI}  \\
		$l=3$ &  0/0  & 1/2 &  2/3  & 3/4 & 4/5  \\
		$l=4$ & 0/0  & 1/2 &  2/3  & 3/4 & 4/5  \\ \hline
	\end{tabular}
	\label{table8}
\end{table}
\begin{center}
\begin{figure}[!ht]
     		\includegraphics[width=6cm, height=4.5cm]{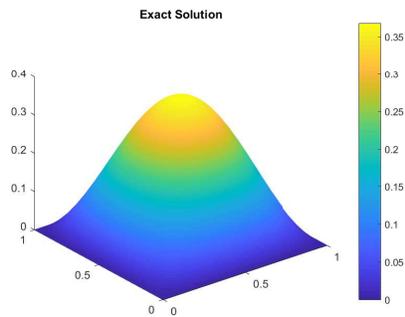}\;\;\centering
      		\caption{Exact solution at $t=1$.}\label{fig5.1.1}
       	\end{figure}
\end{center}

\section{Concluding Remarks}\se
In this paper we have conducted a systematic study for the WG-FEM with local elements $({\cal P}_{k}(K),\;{\cal P}_{j}(\partial K),\;\big[{\cal P}_{l}(K)\big]^2)$.
 For all values of $k, \;j$ and $l$, we have established a theoretical framework for the convergence and error
estimates in the triple bar norm and standard $L^2$ norm. The results are summarized as follows.
\begin{enumerate} [label=(\Alph*)]
	\item For the method of projected element-boundary-discrepancy, we have the following results:
	\begin{enumerate} [label=(\roman*)]
		\item For $l\geq k-1 $ and $j\geq l$, the corresponding WG scheme is stable and the convergence order of $k$ in the triple bar norm and $k+1$ in $L^2$ norm.
		\item For $l\geq k-1$ and $j<l$, the corresponding WG scheme is stable and has convergence order of $s= \min\{k, j\}$ in the triple bar norm and $s+1$ in $L^2$ norm.
		\item For $l<k-1$, the corresponding WG scheme is either unstable or not consistent.
	\end{enumerate}
	\item For the stabilizer with  element-boundary-discrepancy, the following results hold true:
	\begin{enumerate} [label=(\roman*)]
		\item For $l\geq k-1 $, the corresponding WG scheme is stable and the convergence order of $s=\min\{k, j\}$ in the triple bar norm and $s+1$ in $L^2$ norm.
		\item For $l=k-2 $ and $j<k$, the corresponding WG scheme is stable and the convergence order of $s=\min\{k, j\}$ in the triple bar norm. 		
\item For $l=k-2 $ with $j\geq k$ and $l<k-2$ the corresponding WG scheme is solvable but not consistent.
	\end{enumerate}
\end{enumerate}

Numerical results are presented for triangular meshes. It would be very challenging to develop weak Galerkin solver for parabolic equation with weak Galerkin space
 $({\cal P}_{k}(K),\;{\cal P}_{j}(\partial K),\;\big[{\cal P}_{l}(K)\big]^2)$, where $k\ge 1$, $j\ge 0$ and $l\ge 0$ are arbitrary integers, for polygonal meshes. Currently, we are working on it for second order Crank-Nicolson scheme. For the possible extension of this work, we refer to second-order parabolic partial differential equations with non-classic boundary conditions (cf. \cite{dehghan2006computational, dehghan2007one, abbaszadeh2020meshless, dehghan2006numerical, dehghan2005efficient, dehghan2002new}).

\bibliographystyle{abbrv}
\bibliography{WG_Para_Sys_Ref_MMM}

\newpage

\section*{Appendix}
In this section we demonstrate some detailed numerical results for a set of selected values of $k,j$ and $l$.
These result will be in support of the rate of convergence reported in Section \ref{sec7}. The numerical results are organized as follows.
Tables \ref{table9}-\ref{table14} illustrate the table index numbers for the set value of $(k,j,l)$,
and the rest of the tables show the corresponding numerical results. For example, Table \ref{table9} points to
the table index numbers when the stabilizer  ${\cal S}(u_h,v_h)=\sKT h^{-1}_K \big\langle {\cal Q}_m(u_b-u_0|_{\partial K}), {\cal Q}_m(v_b-v_0|_{\partial K})\big\rangle_{\partial K}$ was
employed in the numerical scheme. This table has a fixed value of $k=2$ while $j$ and $l$ are varying. The entry of the table at $(k,j,l) = (2,1,2)$ has value Table \ref{table17},
so that the computational results for $(k,j,l) = (2,1,2) $ should be found in Table \ref{table17}.

\subsection*{Index Tables} \label{In}
The index tables are given in Tables {\ref{table9}}, {\ref{table10}}, {\ref{table11}}, {\ref{table12}}, {\ref{table13}}, and {\ref{table14}}. Here, note that the values in those tables refer to the table number where the computational results are reported.

\begin{table}[!ht]
	\caption{Order of convergence for $k$=2 with stabilizer term
	$~~~~~~~~~~~~~~~~~~~~~~~~ h^{-1}_K \big\langle {\cal Q}_m(u_b-u_0|_{\partial K}),{\cal Q}_m(v_b-v_0|_{\partial K})\big\rangle$ } \centering
\vspace{0.2cm}
	\begin{tabular}{c c c c  c c }
		\hline
		$k=2$ &  $ j=0 $&  $ j=1 $& $j=2$  &$ j=3 $&  $j=4$ \\ \hline
		$l=0$ &     \\
		
		$l=1$ &   Table 15 & Table 16 &    \\

		$l=2$ &    & Table 17 & Table 18   &   \\
			
		$l=3$ &  & & &     Table 19 &         \\
	
		$l=4$ &   &&&&  Table 20  \\ \hline
	\end{tabular}
	\label{table9}
\end{table}
\begin{table}[!ht]
\caption{Order of convergence for $k$=3 with stabilizer term
	$~~~~~~~~~~~~~~~~~~~~~~~~ h^{-1}_K \big\langle {\cal Q}_m(u_b-u_0|_{\partial K}),{\cal Q}_m(v_b-v_0|_{\partial K})\big\rangle$ } \centering
\vspace{0.2cm}
	\begin{tabular}{c c c c  c c }
		\hline
		$k=3$ &  $ j=0 $& $ j=1 $&  $j=2$   &$ j=3 $&   $j=4$ \\ \hline
		$l=0$ &     \\
		
		$l=1$ &   & & &  ~   \\
		
		$l=2$ &   &Table 21 & Table 22   &     \\
		
		$l=3$ &  & & Table 23 &  Table 24       \\

		$l=4$ &   &&&&   Table 25  \\ \hline
	\end{tabular}
	\label{table10}
\end{table}
\begin{table}[!ht]
	\caption{Order of convergence for $k$=4 with stabilizer term
		$~~~~~~~~~~~~~~~~~~~~~~~~ h^{-1}_K \big\langle {\cal Q}_m(u_b-u_0|_{\partial K}),{\cal Q}_m(v_b-v_0|_{\partial K})\big\rangle$ } \centering
	\vspace{0.2cm}
	\begin{tabular}{c c c c  c c }
		\hline
		$k=4$ &  $ j=0 $& $ j=1 $&  $j=2$  &$ j=3 $&  $j=4$ \\ \hline
		$l=0$ &   \\
		
		$l=1$ &   & & & ~ & ~ \\

		$l=2$ &  &  ~ & ~ \\

		$l=3$ &  &  Table 26 & Table 27 & Table 28   \\

		$l=4$ &   &&&   Table 29 &   Table 30  \\ \hline
	\end{tabular}
	\label{table11}
\end{table}
\begin{table}[!ht]
\caption{ $~~~~$ Order of convergence for $k$=2 with stabilizer term\newline
	$~~~~~~~~~~~~~~~~~~~~ h^{-1}_K \big\langle u_b-u_0|_{\partial K},v_b-v_0|_{\partial K}\big\rangle$ } \centering
\vspace{0.2cm}
	\begin{tabular}{c c c c  c c }
		\hline
		$k=2$ &  $ j=0 $& $ j=1 $&  $j=2$ &$ j=3 $&  $j=4$ \\ \hline
		$l=0$ &    \\
		
		$l=1$ &   & Table 31 \\

		$l=2$ &     & & Table 32  &   \\

		$l=3$ &  & & &   Table 33 &         \\

		$l=4$ &   &&&&   Table 34  \\ \hline
	\end{tabular}
	\label{table12}
\end{table}
\begin{table}[!ht]
	\caption{ $~~~~$ Order of convergence for $k$=3 with stabilizer term\newline
		$~~~~~~~~~~~~~~~~~~~~ h^{-1}_K \big\langle u_b-u_0|_{\partial K},v_b-v_0|_{\partial K}\big\rangle$ } \centering
	\vspace{0.2cm}
	\begin{tabular}{c c c c  c c }
		\hline
		$k=3$ &  $ j=0 $& $ j=1 $&  $j=2$  &$ j=3 $&  $j=4$ \\ \hline
		$l=0$ &       \\
		
		$l=1$ &    &  Table 35 & ~ & ~  \\

		$l=2$ &     & & Table 36  &     \\

		$l=3$ &  & & &    Table 37 &   \\

		$l=4$ &   &&&&    Table 38  \\ \hline
	\end{tabular}
	\label{table13}
\end{table}
\begin{table}[!ht]
	\caption{ $~~~~$ Order of convergence for $k$=4 with stabilizer term\newline
		$~~~~~~~~~~~~~~~~~~~~ h^{-1}_K \big\langle u_b-u_0|_{\partial K},v_b-v_0|_{\partial K}\big\rangle$ } \centering
	\vspace{0.2cm}
	\begin{tabular}{c c c c  c c }
		\hline
		$k=4$ &  $ j=0 $& $ j=1 $& $j=2$  &$ j=3 $&  $j=4$ \\ \hline
		$l=0$ &    \\
		$l=1$ & & & ~ & ~ \\
		
		$l=2$ &  & Table 39  & Table 40  & Table 41 & ~   \\
		
		$l=3$ &  & & &   Table 42 & \\
		
		$l=4$ &   &&&&    Table 43 \\ \hline
	\end{tabular}
	\label{table14}
\end{table}

\subsection*{Tables for Computational Results}
All the detailed numerical results are provided in Tables \ref{table15}-\ref{table53}.
No interpretation of the data is necessary as they are virtually self-explanatory. Interested readers are invited to draw their own conclusions from reading these numerical results.
\begin{table}[!ht] 	
\caption{Convergence orders for $k=2,j=0,l=1$ with time step $\tau =10^{-4}$}
	\centering
	\begin{tabular}{c c c c  c c }
		\hline
		$ h $ &  $ \vertiii{e^n} $& $ Order $&  $\|{e^n}\|$  &$ Order $  \\ \hline
		$ 1/4 $ &  8.122260e-01 &   & 9.918705e-02  &  &       \\
		$1/8 $ &  8.458429e-01  & -5.850856e-02& 1.024692e-01  & -4.696613e-02  &   \\
		$ 1/16 $ & 8.547738e-01    & -1.515296e-02 & 1.030492e-01  &  -8.143834e-03  &   \\
		$ 1/32 $ &  8.570361e-01  &  -3.813233e-03 & 1.031780e-01 &  -1.801845e-03  & \\ \hline
	\end{tabular}
	\label{table15}
	\caption{Convergence orders for $k=2,j=1,l=1$ with time step $\tau =10^{-4}$}
	\centering
	\begin{tabular}{c c c c  c c }
		\hline
		$ h $ &  $ \vertiii{e^n} $& $ Order $&  $\|{e^n}\|$  &$ Order $  \\ \hline
		$ 1/4 $ &  7.169166e-02   &   & 6.189540e-03 &  &       \\
		$1/8 $ &  1.805445e-02  & 1.989451e+00 &  7.725189e-04   & 3.002190e+00   &   \\
		$ 1/16 $ & 4.522790e-03   & 1.997070e+00  & 9.652195e-05  &  3.000641e+00  &   \\
		$ 1/32 $ & 1.131375e-03   & 1.999136e+00 & 1.208548e-05  &   2.997582e+00  & \\ \hline
	\end{tabular}
	\label{table16}
\end{table}

\begin{table}[!ht]
	\caption{Convergence orders for $k=2,j=1,l=2$ with time step $\tau =10^{-4}$}
	\centering
	\begin{tabular}{c c c c  c c }
		\hline
		$ h $ &  $ \vertiii{e^n} $& $ Order $&  $\|{e^n}\|$  &$ Order $  \\ \hline
		$ 1/4 $ &  1.652606e-01   &   & 1.071478e-02 &  &       \\
		$1/8 $ &  8.483399e-02   & 9.620281e-01 & 2.906554e-03   & 1.882220e+00   &   \\
		$ 1/16 $ & 4.268569e-02  & 9.908899e-01  & 7.421362e-04   &  1.969554e+00  &   \\
		$ 1/32 $ & 2.137580e-02   & 9.977740e-01 & 1.862606e-04   &   1.994361e+00  & \\ \hline
	\end{tabular}
	\label{table17}
	\caption{Convergence orders for $k=2,j=2,l=2$ with time step $\tau =10^{-5}$}
	\centering
	\begin{tabular}{c c c c  c c }
		\hline
		$ h $ &  $ \vertiii{e^n} $& $ Order $&  $\|{e^n}\|$  &$ Order $  \\ \hline
		$ 1/4 $ &  9.067179e-03   &   & 5.671533e-04 &  &       \\
		$1/8 $   &  1.342686e-03  &  2.755532e+00  & 3.809727e-05   & 3.895979e+00   &   \\
		$ 1/16 $ & 2.412130e-04   &  2.476743e+00  &  2.851774e-06  & 3.739756e+00   &   \\
		$ 1/32 $ & 5.269517e-05   &  2.194565e+00  &  2.672112e-07  &   3.415807e+00 & \\ \hline
	\end{tabular}
	\label{table18}
	\caption{Convergence orders for $k=2,j=3,l=3$ with time step $\tau =10^{-5}$}
	\centering
	\begin{tabular}{c c c c  c c }
		\hline
		$ h $ &  $ \vertiii{e^n} $& $ Order $&  $\|{e^n}\|$  &$ Order $  \\ \hline
		$ 1/4 $ &  5.196970e-03    &   &  1.247593e-04  &  &       \\
		$1/8 $ &  1.269027e-03  & 2.033948e+00 & 1.449317e-05   & 3.105703e+00   &   \\
		$ 1/16 $ & 3.153325e-04   & 2.008777e+00  & 1.767677e-06  &  3.035446e+00   &   \\
		$ 1/32 $ & 7.873437e-05   & 2.001809e+00 & 2.227589e-07   &  2.988299e+00   & \\ \hline
	\end{tabular}
	\label{table19}
	\caption{Convergence orders for $ k=2,j=4,l=4 $ with time step $\tau =10^{-4}$}
	\centering
	\begin{tabular}{c c c c  c c }
		\hline
		$ h $ &  $ \vertiii{e^n} $& $ Order $&  $\|{e^n}\|$  &$ Order $  \\ \hline
		$ 1/4 $ &  5.196970e-03    &   &   6.863976e-04   &  &       \\
		$1/8 $ &  1.269027e-03  & 1.974641e+00 & 8.241034e-05   & 3.046438e+00   &   \\
		$ 1/16 $ & 3.153325e-04   & 1.993420e+00  & 1.012732e-05   &  3.024572e+00   &   \\
		$ 1/32 $ & 7.873437e-05   & 1.998177e+00 & 1.328450e-06    &   2.930437e+00   & \\ \hline
	\end{tabular}
	\label{table20}
	\caption{Convergence orders for $ k=3,j=1,l=2 $ with time step $\tau =10^{-4}$}
	\centering
	\begin{tabular}{c c c c  c c }
		\hline
		$ h $ &  $ \vertiii{e^n} $& $ Order $&  $\|{e^n}\|$  &$ Order $  \\ \hline
		$ 1/4 $ &  1.670987e-01    &   &  1.040284e-02  &  &       \\
		$1/8 $ & 8.570960e-02 & 9.631714e-01 & 2.830643e-03   & 1.877776e+00    &   \\
		$ 1/16 $ & 4.311791e-02 & 9.911696e-01  & 7.233917e-04   &  1.968281e+00  &   \\
		$ 1/32 $ & 2.159118e-02   & 9.978453e-01& 1.815894e-04    &   1.994097e+00    & \\ \hline
	\end{tabular}
	\label{table22}
	\caption{Convergence orders for $ k=3,j=2,l=2 $ with time step $\tau =10^{-5}$}
	\centering
	\begin{tabular}{c c c c  c c }
		\hline
		$ h $ &  $ \vertiii{e^n} $& $ Order $&  $\|{e^n}\|$  &$ Order $  \\ \hline
		$ 1/4 $ &  9.201438e-03    &   &  6.734277e-04  &  &       \\
		$1/8 $ &   1.164020e-03 & 2.982744e+00 & 4.245300e-05   &  3.98758e+00   &   \\
		$ 1/16 $ & 1.459683e-04 & 2.995389e+00  & 2.659047e-06   &  3.986885e+00  &   \\
		$ 1/32 $ & 1.8226353e-05   & 2.998617e+00&  1.733472e-07    &   3.939173e+00    & \\ \hline
	\end{tabular}
	\label{table23}
	\caption{Convergence orders for $ k=3,j=2,l=3 $ with time step $\tau =10^{-4}$}
	\centering
	\begin{tabular}{c c c c  c c }
		\hline
		$ h $ &  $ \vertiii{e^n} $& $ Order $&  $\|{e^n}\|$  &$ Order $  \\ \hline
		$ 1/4 $ &  2.461944e-02   &   &  6.907415e-04  &  &       \\
		$1/8 $ &   6.266297e-03 & 1.974113e+00 & 8.841892e-05   & 2.965719e+00    &   \\
		$ 1/16 $ & 1.572490e-03 & 1.994563e+00  & 1.106445e-05 &  2.998423e+00  &   \\
		$ 1/32 $ & 3.935074e-04  & 1.998588e+00 &  1.452128e-06    &   2.929691e+00   & \\ \hline
	\end{tabular}
	\label{table24}
	\caption{Convergence orders for $ k=3,j=3,l=3 $ with time step $\tau =10^{-6}$}
	\centering
	\begin{tabular}{c c c c  c c }
		\hline
		$ h $ &  $ \vertiii{e^n} $& $ Order $&  $\|{e^n}\|$  &$ Order $  \\ \hline
		$ 1/4 $ & 8.712987e-04    &   & 4.796983e-05   &  &       \\
		$1/8 $ &    6.345320e-05  & 3.779402e+00 &  1.581963e-06  &  4.922340e+00   &   \\
		$ 1/16 $ &  5.576660e-06 & 3.508220e+00 & 5.637939e-08   &   4.810404e+00   &   \\
		$ 1/32 $ & 6.010600e-07  & 3.213821e+00 &  5.462454e-09     & 3.367547e+00   & \\ \hline
	\end{tabular}
	\label{table25}
\end{table}

\begin{table}[!ht]
	\caption{Convergence orders for $ k=3,j=4,l=4 $ with time step $\tau =10^{-6}$}
	\centering
	\begin{tabular}{c c c c  c c }
		\hline
		$ h $ &  $ \vertiii{e^n} $& $ Order $&  $\|{e^n}\|$  &$ Order $  \\ \hline
		$ 1/4 $ &  4.710535e-04   &   & 8.416242e-06  &  &       \\
		$1/8 $ &    5.781054e-05  & 3.026486e+00  & 4.870973e-07    & 4.110894e+00   &   \\
		$ 1/16 $ &  7.207176e-06 &  3.003827e+00  &  2.990930e-08  &  4.025544e+00  &   \\
		$ 1/32 $ &  9.014070e-07 & 2.999184e+00   &  5.234687e-09  &  2.514419e+00   & \\ \hline
	\end{tabular}
	\label{table26}
	\caption{Convergence orders for $ k=4,j=1,l=3 $ with time step $\tau =10^{-4}$}
	\centering
	\begin{tabular}{c c c c  c c }
		\hline
		$ h $ &  $ \vertiii{e^n} $& $ Order $&  $\|{e^n}\|$  &$ Order $  \\ \hline
		$ 1/4 $ &  2.716178e-01   &   &   1.298408e-02   &  &       \\
		$1/8 $ &  1.383182e-01 & 9.735868e-0 & 3.455674e-03   & 1.909705e+00   &   \\
		$ 1/16 $ & 6.946239e-02 & 9.936875e-01  & 8.782836e-04 &  1.976208e+00  &   \\
		$ 1/32 $ & 3.476853e-02 & 9.984499e-01 &  2.202219e-04    &   1.995729e+00    & \\ \hline
	\end{tabular}
	\label{table31}
	\caption{Convergence orders for $ k=4,j=2,l=3 $ with time step $\tau =10^{-4}$}
	\centering
	\begin{tabular}{c c c c  c c }
		\hline
		$ h $ &  $ \vertiii{e^n} $& $ Order $&  $\|{e^n}\|$  &$ Order $  \\ \hline
		$ 1/4 $ &  2.476214e-02   &   &   6.646780e-04  &  &       \\
		$1/8 $ &  6.298864e-03 & 1.974973e+00 & 8.522381e-05   & 2.963327e+00   &   \\
		$ 1/16 $ & 1.580428e-03  & 1.994777e+00  & 1.067046e-05 & 2.997635e+00  &   \\
		$ 1/32 $ &  3.954779e-04  & 1.998646e+00 &  1.405669e-06   &   2.924293e+00    & \\ \hline
	\end{tabular}
	\label{table32}
	\caption{Convergence orders for $ k=4,j=3,l=3 $ with time step $\tau =10^{-6}$}
	\centering
	\begin{tabular}{c c c c  c c }
		\hline
		$ h $ &  $ \vertiii{e^n} $& $ Order $&  $\|{e^n}\|$  &$ Order $  \\ \hline
		$ 1/4 $ & 9.038281e-04   &   &   5.853914e-05  &  &       \\
		$1/8 $ &  5.715505e-05  & 3.983096e+00 & 1.849952e-06  & 4.983842e+00 &   \\
		$ 1/16 $ & 3.583032e-06  & 3.995628e+00 & 5.818479e-08  & 4.990701e+00   &   \\
		$ 1/32 $ &  2.251702e-07  &3.992093e+00 &  5.232165e-09  & 3.475162e+00 & \\ \hline
	\end{tabular}
	\label{table33}
	\caption{Convergence orders for $ k=4,j=3,l=4 $ with time step $\tau =10^{-6}$}
	\centering
	\begin{tabular}{c c c c  c c }
		\hline
		$ h $ &  $ \vertiii{e^n} $& $ Order $&  $\|{e^n}\|$  &$ Order $  \\ \hline
		$ 1/4 $ &  2.705839e-03     &   &   6.050703e-05  &  &       \\
		$1/8 $ &  3.475762e-04 & 2.960675e+00  &  3.841695e-06    & 3.977288e+00  &   \\
		$ 1/16 $ & 4.380892e-05 & 2.988033e+00   & 2.399965e-07  & 4.000658e+000  &   \\
		$ 1/32 $ &  5.491110e-06  & 2.996055e+00  &  1.574981e-08   &   3.929607e+00     & \\ \hline
	\end{tabular}
	\label{table34}
	\caption{Convergence orders for $k=4,j=4,l=4$ with time step $\tau =10^{-6}$}
	\centering
	\begin{tabular}{c c c c c c}
		\hline
		$ h $ &  $\tnorm{e^n}$& $Order$& $\|{e^n}\|$  & $Order$  \\ \hline
		$ 1/4 $ & 7.096378e-05    &  &   3.364898e-06  &  &       \\
		$1/8 $ &  2.503060e-06  & 4.825318e+00  & 7.897632e-08   & 5.412998e+00   &  ~ \\
		$ 1/16 $ & 1.067160e-07   & 4.551845e+00  & 2.586304e-09 & 4.932456e+00  &  ~ \\
		$ 1/32 $ &  2.247061e-08  & 2.247665e+00 &  1.980236e-10   &  3.707147e+00  ~ &~ \\ \hline
	\end{tabular}
	\label{table35}
	\caption{Convergence orders for $ k=2,j=1,l=1 $ with time step $\tau =10^{-4}$}
	\centering
	\begin{tabular}{c c c c  c c }
		\hline
		$ h $ &  $ \vertiii{e^n} $& $ Order $&  $\|{e^n}\|$  &$ Order $  \\ \hline
		$ 1/4 $ &  2.476214e-02   &   &   6.436302e-03   &  &       \\
		$1/8 $ &  6.298864e-03 & 1.438592e+00 & 1.118485e-03    & 2.524686e+00   &   \\
		$ 1/16 $ & 1.580428e-03  & 1.174826e+00  & 2.386911e-04 & 2.228330e+00   &   \\
		$ 1/32 $ &  3.954779e-04  & 1.051661e+00 &  5.674560e-05   &   2.072564e+00    & \\ \hline
	\end{tabular}
	\label{table36}
	\caption{Convergence orders for $ k=2,j=2,l=2 $ with time step $\tau =10^{-5}$}
	\centering
	\begin{tabular}{c c c c  c c }
		\hline
		$ h $ &  $ \vertiii{e^n} $& $ Order $&  $\|{e^n}\|$  &$ Order $  \\ \hline
		$ 1/4 $  &  9.067179e-03   &   &  5.671533e-04  &  &       \\
		$1/8 $   &  1.342686e-03   & 2.755532e+00   & 3.809727e-05  & 3.895979e+00   &   \\
		$ 1/16 $ &  2.412130e-04   & 2.476743e+00   & 2.851774e-06  & 3.739756e+00   &   \\
		$ 1/32 $ &  5.269517e-05   & 2.194565e+00   & 2.672112e-07  & 3.415807e+00   & \\ \hline
	\end{tabular}
	\label{table37}
\end{table}

\begin{table}[!ht]
	\caption{Convergence orders for $ k=2,j=3,l=3 $ with time step $\tau =10^{-5}$}
	\centering
	\begin{tabular}{c c c c  c c }
		\hline
		$ h $ &  $ \vertiii{e^n} $& $ Order $&  $\|{e^n}\|$  &$ Order $  \\ \hline
		$ 1/4 $ &  5.196970e-03   &   &   1.247593e-04  &  &       \\
		$1/8 $ &  1.269027e-03 & 2.033948e+00 & 1.449317e-05  & 3.105703e+00  &   \\
		$ 1/16 $ & 3.153325e-04  & 2.008777e+00  & 1.767677e-06 & 3.035446e+00   &   \\
		$ 1/32 $ &  7.873437e-05  & 2.001809e+00 &  2.227589e-07   &   2.988299e+00    & \\ \hline
	\end{tabular}
	\label{table38}
	\caption{Convergence orders for $ k=2,j=4,l=4 $ with time step $\tau =10^{-4}$}
	\centering
	\begin{tabular}{c c c c  c c }
		\hline
		$ h $ &  $ \vertiii{e^n} $& $ Order $&  $\|{e^n}\|$  &$ Order $  \\ \hline
		$ 1/4 $ &  5.276425e-02    &   &   6.808492e-04    &  &       \\
		$1/8 $ &  1.342498e-02 & 1.974641e+00 & 8.241034e-05     & 3.046438e+00    &   \\
		$ 1/16 $ & 3.371586e-03   & 1.993420e+00  & 1.012732e-05  & 3.024572e+00  &   \\
		$ 1/32 $ &  8.439626e-04   & 1.998177e+00 &  1.328451e-06    &   2.930436e+00    & \\ \hline
	\end{tabular}
	\label{table39}
	\caption{Convergence orders for $ k=3,j=1,l=1 $ with time step $\tau =10^{-4}$}
	\centering
	\begin{tabular}{c c c c  c c }
		\hline
		$ h $ &  $ \vertiii{e^n} $& $ Order $&  $\|{e^n}\|$  &$ Order $  \\ \hline
		$ 1/4 $ &  1.235269e-01     &   &   1.485070e-02   &  &       \\
		$1/8 $ &  4.253941e-02  & 1.537953e+00& 2.061305e-03      & 2.848901e+00    &   \\
		$ 1/16 $ & 1.808968e-02    & 1.233633e+00 & 3.432200e-04  & 2.586353e+00  &   \\
		$ 1/32 $ & 8.601552e-03   & 1.072498e+00 &  7.111218e-05     &  2.270965e+00     & \\ \hline
	\end{tabular}
	\label{table40}
	\caption{Convergence orders for $ k=3,j=2,l=2 $ with time step $\tau =10^{-4}$}
	\centering
	\begin{tabular}{c c c c  c c }
		\hline
		$ h $ &  $ \vertiii{e^n} $& $ Order $&  $\|{e^n}\|$  &$ Order $  \\ \hline
		$ 1/4 $ &  1.048823e-02      &   &    7.276300e-04  &  &       \\
		$1/8 $ &  1.866607e-03   & 2.490281e+00 &  6.014785e-05   & 3.596620e+00  &   \\
		$ 1/16 $ & 4.035579e-04   & 2.209570e+00 & 6.163514e-06  & 3.286688e+00 &   \\
		$ 1/32 $ & 9.652462e-05   & 2.063807e+00 &   8.686333e-07  & 2.826934e+00     & \\ \hline
	\end{tabular}
	\label{table43}
	\caption{Convergence orders for $ k=3,j=3,l=3 $ with time step $\tau =10^{-6}$}
	\centering
	\begin{tabular}{c c c c  c c }
		\hline
		$ h $ &  $ \vertiii{e^n} $& $ Order $&  $\|{e^n}\|$  &$ Order $  \\ \hline
		$ 1/4 $ &  8.712987e-04  &   &    4.796983e-05   &  &       \\
		$1/8 $ &   6.345320e-05 & 3.779402e+00 & 1.581963e-06  &  4.922340e+00 &   \\
		$ 1/16 $ &  5.576661e-06 & 3.508220e+00& 5.637940e-08  & 4.810404e+00   &   \\
		$ 1/32 $ & 6.010604e-07  & 3.213820e+00 & 5.462806e-09 &  3.367454e+00  & \\ \hline
	\end{tabular}
	\label{table44}
	\caption{Convergence orders for $ k=3,j=4,l=4 $ with time step $\tau =10^{-6}$}
	\centering
	\begin{tabular}{c c c c  c c }
		\hline
		$ h $ &  $ \vertiii{e^n} $& $ Order $&  $\|{e^n}\|$  &$ Order $  \\ \hline
		$ 1/4 $ &  4.710535e-04   &   &    8.416242e-06  &  &       \\
		$1/8 $ &  5.781054e-05  & 3.026486e+00 & 4.870973e-07  & 4.110894e+00    &   \\
		$ 1/16 $ & 7.207176e-06   & 3.003827e+00  & 2.990930e-08  &  4.025544e+00  &   \\
		$ 1/32 $ &  9.014071e-07 & 2.999183e+00 & 5.234441e-09  &  2.514487e+00 & \\ \hline
	\end{tabular}
	\label{table45}
	\caption{Convergence orders for $ k=4,j=1,l=2 $ with time step $\tau =10^{-4}$}
	\centering
	\begin{tabular}{c c c c  c c }
		\hline
		$ h $ &  $ \vertiii{e^n} $& $ Order $&  $\|{e^n}\|$  &$ Order $  \\ \hline
		$ 1/4 $ &1.680932e-01 &   &    1.019948e-02  &  &       \\
		$1/8 $ &  8.606377e-02   &  9.657835e-01  & 2.813409e-03 & 1.858104e+00  &   \\
		$ 1/16 $ & 4.329160e-02   & 9.913190e-01  & 7.223234e-04  & 1.961603e+00 &   \\
		$ 1/32 $ & 2.167803e-02   &  9.978530e-01 &  1.815416e-04 & 1.992344e+00  & \\ \hline
	\end{tabular}
	\label{table48}
	\caption{Convergence orders for $ k=4,j=2,l=2 $ with time step $\tau =10^{-4}$ }
	\centering
	\begin{tabular}{c c c c  c c }
		\hline
		$ h $ &  $ \vertiii{e^n} $& $ Order $&  $\|{e^n}\|$  &$ Order $  \\ \hline
		$ 1/4 $ & 1.659611e-02 &   & 1.913745e-03   &  &       \\
		$1/8 $ & 2.746179e-03 &  2.595348e+00 &  1.336710e-04  & 3.839640e+00 &   \\
		$ 1/16 $ &  5.644401e-04  & 2.282533e+00  & 1.117750e-05  & 3.580018e+00  &   \\
		$ 1/32 $ & 1.323725e-04  & 2.092217e+00  &  1.258791e-06  & 3.150486e+00 & \\ \hline
	\end{tabular}
	\label{table49}
\end{table}

\begin{table}[t]
	\caption{Convergence orders for $ k=4,j=3,l=2 $ with time step $\tau =10^{-4}$}
	\centering
	\begin{tabular}{c c c c  c c }
		\hline
		$ h $ &  $ \vertiii{e^n} $& $ Order $&  $\|{e^n}\|$  &$ Order $  \\ \hline
		$ 1/4 $ &  1.475717e-02  &   &  1.880159e-03   &  &       \\
		$1/8 $ &   1.839962e-03  & 3.003668e+00  & 1.154356e-04   & 4.025694e+00   &   \\
		$ 1/16 $ & 2.298950e-04    & 3.000629e+00  & 7.199269e-06  & 4.003094e+00  &   \\
		$ 1/32 $ &  2.881565e-05   & 2.996050e+00 &  6.649049e-07 &  3.436630e+00 & \\ \hline
	\end{tabular}
	\label{table50}
	\caption{Convergence orders for $ k=4,j=3,l=3 $ with time step $\tau =10^{-6}$ }
	\centering
	\begin{tabular}{c c c c  c c }
		\hline
		$ h $ &  $ \vertiii{e^n} $& $ Order $&  $\|{e^n}\|$  &$ Order $  \\ \hline
		$ 1/4 $ & 1.036100e-03  &   &    6.412737e-05  &  &       \\
		$1/8 $ &  9.113811e-05   & 3.506965e+00  &   2.618642e-06        & 4.614049e+00   &   \\
		$ 1/16 $ &   9.779329e-06      & 3.220247e+00  & 1.323056e-07 & 4.306873e+00  &   \\
		$ 1/32 $ & 1.166344e-06    & 3.067742e+00  &  9.121470e-09 & 3.858464e+00 & \\ \hline
	\end{tabular}
	\label{table52}
	\caption{Convergence orders for $ k=4,j=4,l=4 $ with time step $\tau =10^{-6}$}
	\centering
	\begin{tabular}{c c c c  c c }
		\hline
		$ h $ &  $ \vertiii{e^n} $& $ Order $&  $\|{e^n}\|$  &$ Order $  \\ \hline
		$ 1/4 $ & 7.096379e-05    &   &   3.3648980e-06  &  &       \\
		$1/8 $ &  2.503069e-06  & 4.825313e+00  & 7.895302e-08   & 5.412998e+00   &   \\
		$ 1/16 $ & 1.067212e-07   & 4.551780e+00  & 2.582612e-09 & 4.934091e+00  &   \\
		$ 1/32 $ & 2.247073e-08   & 2.247728e+00&  1.988624e-10   &  3.698988e+00    & \\ \hline
	\end{tabular}
	\label{table53}
\end{table}

\end{document}